\documentclass[11pt]{article}

\usepackage[letterpaper,margin=1in]{geometry}

\usepackage[T1]{fontenc}
\usepackage[utf8]{inputenc}
\usepackage{lmodern}

\usepackage{amsmath}
\usepackage{amssymb}
\usepackage{amsfonts}
\usepackage{amsthm}
\usepackage{amsopn}
\usepackage{mathtools}
\usepackage{mathrsfs}
\usepackage{bm}
\usepackage{stmaryrd}
\usepackage{braket}
\usepackage{autobreak}
\allowdisplaybreaks

\usepackage{array}
\usepackage{multirow}
\usepackage{makecell}
\usepackage{diagbox}
\usepackage{graphicx}
\usepackage{epstopdf}
\usepackage{float}
\usepackage{subfigure}     
\usepackage{pgfplots}
\usepackage{capt-of}
\pgfplotsset{compat=newest}

\usepackage{algorithmic}

\usepackage[titletoc,title]{appendix}
\usepackage{xspace}
\usepackage{etoolbox}

\usepackage[hidelinks]{hyperref}

\usepackage{cleveref}

\newcommand{\vertiii}[1]{{\left\vert\kern-0.25ex\left\vert\kern-0.25ex\left\vert #1
        \right\vert\kern-0.25ex\right\vert\kern-0.25ex\right\vert}}

\theoremstyle{plain}
\newtheorem{theorem}{Theorem}[section]
\newtheorem{lemma}[theorem]{Lemma}

\newtheorem{corollary}[theorem]{Corollary}

\theoremstyle{definition}

\theoremstyle{remark}
\newtheorem{remark}[theorem]{Remark}
\newtheorem{assumption}[theorem]{Assumption}

\crefname{assumption}{Assumption}{Assumptions}
\Crefname{assumption}{Assumption}{Assumptions}
\crefname{hypothesis}{Hypothesis}{Hypotheses}
\Crefname{hypothesis}{Hypothesis}{Hypotheses}

\hypersetup{
  pdftitle={Sharp L2 error estimates for the stabilized nonsymmetric Nitsche method on convex polytopes},
  pdfauthor={Gang Chen, Chaoran Liu, Yangwen Zhang}
}

\begin{document}

\makeatletter
\providecommand\firstaid@cref@smugglelabel{\let\cref@currentlabel\cref@gcurrentlabel@temp}
\providecommand\firstaid@cref@updatelabeldata[1]{%
  \cref@constructprefix{#1}{\cref@result}%
  \@ifundefined{cref@#1@alias}{\def\@tempa{#1}}{\def\@tempa{\csname cref@#1@alias\endcsname}}%
  \protected@xdef\cref@gcurrentlabel@temp{%
    [\@tempa][\arabic{#1}][\cref@result]%
    \csname p@#1\endcsname\csname the#1\endcsname}%
  \aftergroup\firstaid@cref@smugglelabel
}
\AddToHook{label}[firstaid/cleveref]{%
  \ifx\@currentcounter\@empty\else
    \iftag@\else
      \firstaid@cref@updatelabeldata{\@currentcounter}%
    \fi
  \fi
}
\makeatother

\title{Sharp and unified $L^2$ error estimates for the
       nonsymmetric Nitsche method on convex polytopes%
\thanks{Gang Chen and Chaoran Liu are supported by the National Natural
        Science Foundation of China (NSFC) under Grant No.\ 121713413 and
        12422115, and by the Jiangsu Provincial Scientific Research Center
        of Applied Mathematics under Grant No.\ BK20233002.
        Yangwen Zhang is supported by the National Science Foundation
        DMS-2111315.}}

\author{%
  Gang Chen\thanks{Department of Mathematics, Sichuan University
    (\texttt{cglwdm@scu.cn}).}
  \and
  Chaoran Liu\thanks{Department of Mathematics, Sichuan University
    (\texttt{mathliuchr@stu.scu.edu.cn}).}
  \and
  Yangwen Zhang\thanks{Department of Mathematics, University of Louisiana
    at Lafayette, Lafayette, LA
    (\texttt{yangwen.zhang@louisiana.edu}).}%
}

\date{\today}

\maketitle

\begin{abstract}
Nitsche's method weakly imposes Dirichlet boundary conditions, but its nonsymmetric variant has long shown a gap between theory and computation: the classical $L^2$~analysis under $H^{k+1}$~regularity predicts a half-order convergence loss, whereas numerical experiments on smooth test problems consistently produce the optimal rate.  Whether this discrepancy reflects a limitation of the analysis or an essential feature of the method has remained an open question.

On bounded convex polytopes in two and three dimensions, we prove a unified, regularity-dependent $L^2$~error estimate valid across the entire penalty scale $h^{-\alpha}$:
\begin{align*}
	\|u-u_h\|_{L^2(\Omega)}\le C h^r |u|_{W^{k+1,p}(\Omega)},\qquad
	r=\min\bigl\{k+1,\,k+\max\{1,\alpha\}-1/p\bigr\}. 
\end{align*}

Numerical experiments in two and three dimensions, on a one-parameter family of manufactured solutions with tunable regularity, demonstrate the sharpness of the estimate and resolve the open question.  First, under merely $H^{k+1}$~regularity the half-order loss is essential; second, the optimal convergence consistently observed on smooth test problems is therefore explained by their full $W^{k+1,\infty}$~regularity, not by a limitation of the standard analysis.  The theory identifies $\alpha\ge 1+1/p$ or $p=\infty$ as the sharp threshold for recovering the optimal rate $h^{k+1}$, and the experiments confirm this if-and-only-if condition in both dimensions.
\end{abstract}


\medskip

	\section{Introduction}
\label{sec:introduction}

Let $\Omega\subset\mathbb R^d$, $d\in\{2,3\}$, be a bounded convex polygonal (resp.\ polyhedral) domain with boundary $\Gamma:=\partial\Omega$ and unit outward normal $n$. We consider the Poisson problem
\begin{equation}\label{org}
	-\Delta u = f \text{ in }\Omega, \qquad u = g \text{ on }\Gamma. 
\end{equation}

Nitsche's method, introduced in~\cite{MR341903}, and now commonly viewed as an alternative to penalty approaches such as Babu\v{s}ka's penalty method~\cite{MR0351118}, provides a consistent variational framework for imposing the Dirichlet condition weakly, by augmenting the bilinear form with boundary terms that enforce $u=g$ without requiring the trial space to satisfy it. This is convenient when the mesh does not conform to the physical boundary or to material interfaces, as in unfitted, fictitious domain, cut, and interface methods~\cite{MR1941489}. The method has since become a standard tool across applications, including incompressible flows with slip boundary conditions~\cite{MR4379970}, fluid--structure interaction~\cite{MR4383066}, virtual element and hybridizable discontinuous Galerkin methods~\cite{MR4732158,MR4752199}, elasticity and contact problems~\cite{MR3719025}, isogeometric multipatch coupling~\cite{MR4449551}, and immersed, cut, and extended finite element methods~\cite{MR4668240}. 

Let $\{\mathcal T_h\}$ be a shape-regular, quasi-uniform family of conforming triangulations of $\Omega$ with mesh size $h:=\max_{T\in\mathcal T_h}\operatorname{diam}(T)$.  For a fixed integer $k\ge 1$, let $\mathcal P_k(T)$ denote the space of polynomials of total degree at most $k$ on $T$, and set
\begin{align*}
	V_h := \{v_h\in H^1(\Omega): v_h|_T\in \mathcal P_k(T) \text{ for all } T\in\mathcal T_h\}.
\end{align*}
Writing $\partial_n v := \nabla v\cdot n$ for the outward normal derivative on $\Gamma$, the Nitsche method with parameters $\beta$ and $\gamma$ seeks $u_h\in V_h$ such that
\begin{align}
	\begin{split}
		(\nabla u_h,\nabla v_h)_\Omega
		&-\langle \partial_n u_h,v_h\rangle_\Gamma
		+\beta\langle u_h,\partial_n v_h\rangle_\Gamma
		+\gamma\langle u_h,v_h\rangle_\Gamma \\
		&= (f,v_h)_\Omega +\beta\langle g,\partial_n v_h\rangle_\Gamma +\gamma\langle g,v_h\rangle_\Gamma \quad \forall v_h\in V_h.
	\end{split}
\end{align}
Here,  $(\cdot,\cdot)_\Omega$ and $\langle\cdot,\cdot\rangle_\Gamma$ denote the $L^2$ inner products on $\Omega$ and $\Gamma$, respectively.

The choice $\beta=-1$ gives the symmetric variant, which with $\gamma=c_0h^{-1}$ for $c_0$ sufficiently large is coercive and admits the optimal $L^2$ estimate $\|u-u_h\|_{L^2(\Omega)}\le Ch^{k+1}|u|_{H^{k+1}(\Omega)}$ via Aubin--Nitsche duality.  The present paper concerns the \emph{nonsymmetric} variant $\beta=+1$, studied across the entire penalty scale $\gamma=h^{-\alpha}$ with $\alpha\in\{-\infty\}\cup\mathbb R$, where $\alpha=-\infty$ denotes the penalty-free endpoint.  This scale unifies three sub-regimes traditionally treated separately:
Freund and Stenberg's classical stabilized scaling $\alpha=1$ \cite{Freund1995OnWI},
the strong-penalty range $\alpha>1$, whose subrange $\alpha\ge 2$
restores the optimal Aubin--Nitsche rate in the standard
$H^{k+1}$ analysis, and Burman's penalty-free endpoint
$\alpha=-\infty$~\cite{MR3022206}.  The nonsymmetric variant is attractive because it is less sensitive to the penalty size than the symmetric form, and remains stable in regimes where the symmetric form is not.

The energy-norm theory is settled, but the $L^2$ theory is not.  Freund and Stenberg~\cite[Theorem~2]{Freund1995OnWI} proved the optimal energy-norm estimate
\begin{align*}
	\|\nabla(u-u_h)\|_{L^2(\Omega)}\le Ch^k|u|_{H^{k+1}(\Omega)},
\end{align*}
while the standard duality argument gives only the half-order $L^2$ rate
\begin{align*}
	\|u-u_h\|_{L^2(\Omega)}\le Ch^{k+1/2}|u|_{H^{k+1}(\Omega)}. 
\end{align*}
Burman~\cite{MR3022206} obtained the same half-order loss for the penalty-free nonsymmetric method.  Yet numerical experiments on smooth test problems consistently produce the full $h^{k+1}$ rate~\cite{MR4744101,MR4994515}, and in his 2012 paper Burman claimed:
\begin{quote}\itshape ``We have not managed to construct an example exhibiting the suboptimal convergence order of the Nitsche method.''\upshape~\cite[\S 8.1, p.~1977]{MR3022206}\end{quote}
Whether the half-order loss reflects a limitation of the analysis or an essential feature of the method has remained an open question.

We resolve this question by proving a unified $L^2$ error estimate (\Cref{thm:unified-L2}) that, for $u\in W^{k+1,p}(\Omega)$ with $p\in[2,\infty]$ and every $\alpha\in\{-\infty\}\cup\mathbb R$, takes the form
\begin{equation*}
	\|u-u_h\|_{L^2(\Omega)} \le C h^r |u|_{W^{k+1,p}(\Omega)},\quad  
	r= \min\bigl\{k+1,\,k+\max\{1,\alpha\}-1/p\bigr\},
\end{equation*}
with $C$ independent of $h$. In the weak-penalty range $\alpha<1$ (including Burman's penalty-free endpoint $\alpha=-\infty$), the exponent reduces to $r=k+1-1/p$ independently of the particular $\alpha$, recovering Burman's half-order rate at $p=2$ and improving continuously to the optimal rate $k+1$ at $p=\infty$.  For the classical stabilized point $\alpha=1$, the same exponent $k+1-1/p$ holds for all $p\in[2,\infty]$.  In the strong-penalty range $\alpha>1$, the rate becomes $k+\alpha-1/p$ in the strong-suboptimal subregion $1<\alpha<1+1/p$ and saturates at the optimal $k+1$ once $\alpha\ge 1+1/p$ or $p=\infty$.

Numerical experiments in \Cref{sec:numerics} on a parametric family of manufactured solutions with tunable regularity, in both $d=2$ and $d=3$, realize $h^{k+1/2}$ at both $\alpha=1$ and $\alpha=-\infty$, establishing sharpness at the $H^{k+1}$ endpoint; at the smooth endpoint $W^{k+1,\infty}$ the same experiments recover the optimal rate $h^{k+1}$ for every $\alpha$, confirming that the half-order loss is intrinsic to the nonsymmetric method at the $H^{k+1}$ endpoint. A separate experiment confirms the sharp if-and-only-if threshold $\alpha\ge 1+1/p$ (or $p=\infty$) for the optimal rate. Together these results resolve the open question: the half-order loss is essential under merely $H^{k+1}$ regularity, while the optimal rate observed on smooth test problems is explained by their full $W^{k+1,\infty}$ regularity rather than by any limitation of the standard analysis.

The half-order phenomenon is connected to a general difficulty: boundary terms may reduce the order predicted by standard duality arguments unless the boundary trace error is controlled at the optimal rate. For adjoint-inconsistent interior penalty discontinuous Galerkin methods, the loss of $L^2$ accuracy is genuine in both theory and computation~\cite{MR1702201,MR1885715}. Related half-order phenomena occur in the analysis of Dirichlet boundary control problems, where conforming, mixed, and HDG analyses often lose half an order for smoother data even though computations show the full rate~\cite{MR2558321,MR3070527,MR4167066,HuMateosSinglerYangwenZhang2018,MR3831243,MR3992054}. Chen et al.~\cite{MR4230430} used maximum norm HDG estimates to recover the missing half order under $W^{k+1,\infty}$ regularity. The same mechanism underlies \Cref{thm:main-L2}: maximum norm information closes the gap.

\section{The traditional \texorpdfstring{$H^{k+1}$}{Hk+1} approach and its obstruction}\label{sec:preliminaries}

\begin{assumption}\label{assump:standing}
	Throughout the paper:
	\begin{enumerate}
		\item[(i)] $\Omega\subset\mathbb R^d$, $d\in\{2,3\}$, is a bounded convex polygonal domain when $d=2$ and a bounded convex polyhedral domain when $d=3$. 
		\item[(ii)] $\{\mathcal T_h\}$ is a conforming, shape regular, quasi uniform family of triangulations of $\Omega$, with shape regularity constant $\sigma$, element diameters $h_T:=\operatorname{diam}(T)$, and $h:=\max_{T\in\mathcal T_h}h_T$.
		\item[(iii)] $k\ge1$ is fixed.  The generic constant $C>0$ is independent of $h$ and of the data; it may depend on $\Omega$, $d$, $k$, $\sigma$, the fixed exponent $\alpha$, and, when $W^{k+1,p}$ estimates are used, on the fixed index $p\in[2,\infty]$.
	\end{enumerate}
\end{assumption}

For an open set $U\subset\mathbb R^d$ or a piecewise smooth subset $U\subset\Gamma$, we use the standard Sobolev spaces $W^{m,p}(U)$, $1\le p\le\infty$, and write $H^m(U)=W^{m,2}(U)$. Norms and seminorms are denoted by $\|\cdot\|_{W^{m,p}(U)}$, $|\cdot|_{W^{m,p}(U)}$, $\|\cdot\|_{H^m(U)}$, and $|\cdot|_{H^m(U)}$, with the convention $\|\cdot\|_{L^p(U)}:=\|\cdot\|_{W^{0,p}(U)}$. For an integer $m\ge 0$, $D^m u$ denotes the vector of all $m$th order weak partial derivatives of $u$, with multi-indices listed once each with multiplicity, and $|D^m u|$ its Euclidean (Frobenius) norm. The subspace of $W^{1,p}(U)$ of functions with vanishing trace on $\partial U$ is denoted $W^{1,p}_0(U)$; in particular, $H_0^1(U):=W^{1,2}_0(U)$.  The fractional Sobolev space $H^{1/2}(\Gamma)$ is the standard trace space of
$H^1(\Omega)$ on the Lipschitz boundary $\Gamma$ [24, \S1.3.3]. We use the
standard surface Sobolev spaces $H^s(\Gamma)$, $0 \le s \le 1$. Since $\Gamma$
is a finite union of flat faces,
\[
\Gamma = \bigcup_{i=1}^{N} \overline{F_i},
\]
the surface gradient is defined facewise. On each face $F_i$, fix the
unit outward normal $n_i$, and for $w \in H^1(F_i)$ let $\nabla_{F_i} w$ denote
the intrinsic weak gradient of $w$ along $F_i$, a tangential vector field on
$F_i$. We set
\[
(\nabla_\Gamma v)|_{F_i} := \nabla_{F_i}(v|_{F_i}).
\]
If $v$ is the trace of a function $\widetilde v$ smooth up to $F_i$, this agrees
with the tangential projection of the ambient gradient,
\[
\nabla_\Gamma v|_{F_i} = \nabla \widetilde v|_{F_i}
- \bigl(n_i \cdot \nabla \widetilde v|_{F_i}\bigr) n_i.
\]
The surface Sobolev space is
\[
H^1(\Gamma) := \bigl\{ v \in L^2(\Gamma) : v|_{F_i} \in H^1(F_i),\ \text{with matching traces on common interfaces} \bigr\},
\]
where the common interfaces are vertices when $d=2$ and edges when $d=3$. Its
norm is equivalent to
\[
\|v\|_{H^1(\Gamma)}^2 \simeq \|v\|_{L^2(\Gamma)}^2
+ \sum_{i=1}^{N} \|\nabla_{F_i}(v|_{F_i})\|_{L^2(F_i)}^2.
\]
The intermediate space is obtained by interpolation,
\[
H^{1/2}(\Gamma) = \bigl[L^2(\Gamma), H^1(\Gamma)\bigr]_{1/2}.
\]

\subsection{Mesh, finite element space, and the nonsymmetric form}\label{sec:mesh-fe-spaces}

Let $\mathcal E_h^\Gamma$ be the set of boundary faces of $\mathcal T_h$; these are edges in two dimensions and triangular faces in three dimensions.  If $E\in\mathcal E_h^\Gamma$, then $T_E$ denotes the unique element having $E$ as a boundary face.  The conforming Lagrange finite element space is
\begin{align*}
	V_h:=\{v_h\in H^1(\Omega):v_h|_T\in\mathcal P_k(T)\text{ for all }T\in\mathcal T_h\},
\end{align*}
with homogeneous subspace and boundary trace space
\begin{align*}
	V_h^0:=V_h\cap H_0^1(\Omega),\qquad M_h:=V_h|_\Gamma. 
\end{align*}
Thus $M_h$ is a continuous piecewise polynomial space of degree $\le k$ on the boundary mesh.  For $v_h\in V_h$, $\partial_n v_h$ on $E\in\mathcal E_h^\Gamma$ is understood elementwise, namely 
$$\partial_n v_h|_E:=\nabla(v_h|_{T_E})\cdot n.$$

We now introduce a one-parameter family of nonsymmetric Nitsche formulations that includes, as special cases, both the classical stabilized scheme and Burman's penalty-free method. The classical method of Nitsche~\cite{MR341903} corresponds to the penalty coefficient $\gamma=c_0h^{-1}$ with $c_0>0$ a fixed constant; at the opposite extreme, Burman~\cite{MR3022206} analyzed the unstabilized case $\gamma=0$. To unify the analysis we consider the entire scale $\gamma=c_0h^{-\alpha}$, with $\alpha\in\{-\infty\}\cup\mathbb R$, where $\alpha=-\infty$ denotes the penalty-free endpoint and the boundary penalty term is understood to vanish at this value. Since a fixed positive coefficient $c_0$ in front of $h^{-\alpha}$ only changes constants, it is normalized to one. The nonsymmetric bilinear form on this scale is
\begin{equation}\label{eq:ah-def}
	a(w,v):=(\nabla w,\nabla v)_\Omega-\langle\partial_n w,v\rangle_\Gamma
	+\langle w,\partial_n v\rangle_\Gamma+h^{-\alpha}\langle w,v\rangle_\Gamma,
\end{equation}
and the nonsymmetric Nitsche approximation seeks $u_h\in V_h$ such that
\begin{equation}\label{eq:prelim-nitsche}
	a(u_h,v_h)=(f,v_h)_\Omega+\langle g,\partial_n v_h\rangle_\Gamma+h^{-\alpha}\langle g,v_h\rangle_\Gamma
	\qquad\forall v_h\in V_h. 
\end{equation}
For the exact solution $u$ of~\eqref{org}, consistency gives the same identity with $u_h$ replaced by $u$; hence Galerkin orthogonality holds:
\begin{equation}\label{eq:galerkin-orthogonality}
	a(u-u_h,v_h)=0\qquad\forall v_h\in V_h. 
\end{equation}
The natural stability norm is the reference Nitsche norm
\begin{equation}\label{eq:reference-nitsche-norm}
	\vertiii{v}_h^2:=\|\nabla v\|_{L^2(\Omega)}^2+h^{-\alpha_\ast}\|v\|_{L^2(\Gamma)}^2,\qquad \alpha_\ast:=\max\{1,\alpha\}.
\end{equation}

The two nonsymmetric boundary terms cancel on the diagonal, so for $\alpha\ge1$,
\begin{equation}\label{eq:diag-cancel}
	a(v_h,v_h)=\|\nabla v_h\|_{L^2(\Omega)}^2+h^{-\alpha}\|v_h\|_{L^2(\Gamma)}^2
	=\vertiii{v_h}_h^2\qquad\forall v_h\in V_h. 
\end{equation}
For $\alpha<1$, the needed weak-penalty inf--sup estimate is stated next. The proof is given in \Cref{app:weak-infsup}.

\begin{theorem}\label{thm:weak-infsup}
	Assume Assumption~\ref{assump:standing}. Let $\alpha\in\{-\infty\}\cup(-\infty,1)$.  Then there exist $h_0>0$ and a constant $c>0$, independent of $h$, of $\alpha$, and of $v_h$, such that, for every mesh with $h<h_0$ and every $v_h\in V_h$, 
	\begin{equation}\label{eq:weak-infsup}
		c\,\vertiii{v_h}_h
		\le
		\sup_{0\ne w_h\in V_h}\frac{|a(v_h,w_h)|}{\vertiii{w_h}_h}.
	\end{equation}
\end{theorem}

\subsection{Basic \texorpdfstring{$H^{k+1}$}{Hk+1} interpolation and projection estimates}

We collect the estimates needed for the traditional argument.  The local scaled trace inequality~\cite[Lemma~12.15]{Ern2021} states that, for every element $T$ and face $E\subset\partial T$,
\begin{equation}\label{eq:scaled-trace}
	\|v\|_{L^2(E)}^2\le C\bigl(h_T^{-1}\|v\|_{L^2(T)}^2+h_T\|\nabla v\|_{L^2(T)}^2\bigr),\qquad v\in H^1(T).
\end{equation}
Together with the inverse inequality~\cite[Lemma~4.5.3, Eq.~(4.5.4), p.~111]{Brenner2008}, this gives
\begin{equation}\label{eq:grad-inverse-trace}
	\|\partial_n v_h\|_{L^2(E)}\le C h_T^{-1/2}\|\nabla v_h\|_{L^2(T)},\qquad v_h\in\mathcal P_k(T). 
\end{equation}
Throughout the paper we shall repeatedly use the following elliptic regularity for the duality argument on the convex polytope $\Omega$~\cite[Theorem~3.2.1.2, p.~147]{MR3396210},
\begin{equation}\label{dual}
	-\Delta z=\phi\text{ in }\Omega,\qquad z=0\text{ on }\Gamma
	\quad\Longrightarrow\quad
	\|z\|_{H^2(\Omega)}\le C\|\phi\|_{L^2(\Omega)}. 
\end{equation}
The normal-trace operator $\partial_n:H^2(\Omega)\cap H_0^1(\Omega)\to H^{1/2}(\Gamma)$ is bounded in the Lions--Magenes/Grisvard sense~\cite[Chapter~1, Theorem~9.4]{Lions1972},~\cite[\S\S~1.3.3 and~1.5, in particular Theorem~1.5.2.1]{MR3396210}, so that
\begin{equation}\label{trace_theorem}
	\|\partial_n z\|_{L^2(\Gamma)}\le\|\partial_n z\|_{H^{1/2}(\Gamma)}\le C\|z\|_{H^2(\Omega)},\qquad z\in H^2(\Omega)\cap H_0^1(\Omega). 
\end{equation}

The standard Lagrange interpolant $I_h$, well-defined on $H^s(\Omega)$ for $2\le s\le k+1$ since $H^s\hookrightarrow C^0(\overline\Omega)$ when $d\le 3$, satisfies the element-wise interpolation estimate~\cite[Theorem~4.4.4, p.~105]{Brenner2008}
\begin{equation}\label{eq:H-interp-volume}
	|v-I_hv|_{H^m(T)}\le C h_T^{s-m}|v|_{H^s(T)},\qquad 0\le m\le s,\ T\in\mathcal T_h,\ v\in H^s(\Omega). 
\end{equation}

We shall use the global $L^2(\Gamma)$ projection onto the continuous trace space $M_h$:
\[
\Pi_h^\Gamma:L^2(\Gamma)\to M_h,
\qquad
\langle \varphi-\Pi_h^\Gamma\varphi,\mu_h\rangle_\Gamma=0\quad\forall\mu_h\in M_h. 
\]
It satisfies the $L^2(\Gamma)$ stability
\begin{equation}\label{eq:Pi-L2-stab}
	\|\Pi_h^\Gamma\varphi\|_{L^2(\Gamma)}\le\|\varphi\|_{L^2(\Gamma)},\qquad\varphi\in L^2(\Gamma),
\end{equation}
and the approximation property~\cite[Chapter~10]{Steinbach2008}
\begin{equation}\label{eq:Pi-H12-approx}
	\|\varphi-\Pi_h^\Gamma\varphi\|_{L^2(\Gamma)}\le Ch^{s}\|\varphi\|_{H^{s}(\Gamma)},\qquad\varphi\in H^{s}(\Gamma),\ {0\le s\le 1}.
\end{equation}
Combining~\eqref{eq:Pi-L2-stab}, \eqref{eq:Pi-H12-approx} with $s=1/2$, and the trace bound~\eqref{trace_theorem}, we obtain
\begin{equation}\label{eq:normal-Pi-est}
	\|\Pi_h^\Gamma\partial_n z\|_{L^2(\Gamma)}
	+h^{-1/2}\|\partial_n z-\Pi_h^\Gamma\partial_n z\|_{L^2(\Gamma)}
	\le C\|z\|_{H^2(\Omega)},\quad z\in H^2(\Omega)\cap H_0^1(\Omega). 
\end{equation}

We shall also use the boundary patch $\mathcal T_h^\Gamma:=\{T\in\mathcal T_h:\overline T\cap\Gamma\ne\emptyset\}$ and the corresponding boundary strip
\begin{equation}\label{eq:strip-def}
	\Omega_h^\Gamma:=\operatorname{int}\bigcup_{T\in\mathcal T_h^\Gamma}\overline T,
	\qquad |\Omega_h^\Gamma|\le Ch,
\end{equation}
and the following boundary-supported discrete lifting on it. {The proof is given in \Cref{app:lem:disc-lift}.}

\begin{lemma}\label{lem:disc-lift}
	Under Assumption~\ref{assump:standing}, there exists a linear operator $E_h:M_h\to V_h$ such that
	\begin{equation}\label{eq:disc-lift-trace}
		(E_h\mu_h)|_\Gamma=\mu_h,
		\qquad
		\operatorname{supp}(E_h\mu_h)\subset\overline{\Omega_h^\Gamma},
		\qquad
		E_h\mu_h=0\text{ a.e.\ in }\Omega\setminus\Omega_h^\Gamma,
	\end{equation}
	and
	\begin{align}
		\|\nabla E_h\mu_h\|_{L^2(\Omega)}&\le C h^{-1/2}\|\mu_h\|_{L^2(\Gamma)}, \label{eq:disc-lift-W12}\\
		\|\partial_n E_h\mu_h\|_{L^2(\Gamma)}&\le C h^{-1}\|\mu_h\|_{L^2(\Gamma)}.  \label{eq:disc-lift-normal}
	\end{align}
\end{lemma}

\begin{corollary}\label{cor:dual-lift-scales}
	For any $z\in H^2(\Omega)\cap H_0^1(\Omega)$, set $v_h:=E_h(\Pi_h^\Gamma\partial_nz)\in V_h$.  Then
	\begin{align}
		\|\Pi_h^\Gamma\partial_nz\|_{L^2(\Gamma)} &\le C\|z\|_{H^2(\Omega)}, \label{eq:Pi-trace-bound}\\
		\|\nabla v_h\|_{L^2(\Omega)} &\le Ch^{-1/2}\|z\|_{H^2(\Omega)}, \label{eq:dual-lift-W12}\\
		\|\partial_n v_h\|_{L^2(\Gamma)} &\le Ch^{-1}\|z\|_{H^2(\Omega)}.  \label{eq:dual-lift-normal}
	\end{align}
\end{corollary}

\begin{proof}
	The $L^2(\Gamma)$ stability~\eqref{eq:Pi-L2-stab} and the trace bound~\eqref{trace_theorem} give~\eqref{eq:Pi-trace-bound}.  Combining~\eqref{eq:Pi-trace-bound} with~\eqref{eq:disc-lift-W12} gives~\eqref{eq:dual-lift-W12}; combining~\eqref{eq:Pi-trace-bound} with~\eqref{eq:disc-lift-normal} gives~\eqref{eq:dual-lift-normal}.
\end{proof}

For $v\in H^1(\Omega)$ with $v|_\Gamma\in L^2(\Gamma)$, we shall use the Ritz projection $\mathcal R_hv\in V_h$ with projected boundary trace, defined by
\begin{equation}\label{eq:nh-ritz-def}
	(\mathcal R_hv)|_\Gamma=\Pi_h^\Gamma(v|_\Gamma),
	\qquad
	(\nabla(v-\mathcal R_hv),\nabla w_h)_\Omega=0\quad\forall w_h\in V_h^0. 
\end{equation}
Then
\begin{equation}\label{eq:Ritz-boundary-orthogonality}
	\langle v-\mathcal R_hv,\mu_h\rangle_\Gamma=0\qquad\forall\mu_h\in M_h. 
\end{equation}
This boundary orthogonality is the reason for using the global projection in~\eqref{eq:nh-ritz-def}.  The following lemma provides the $L^2$ error estimates for the Ritz projection and its proof is given in \Cref{app:ritz_pro_L2err}.

\begin{lemma}\label{ritz_pro_L2err}
	For $v\in H^{k+1}(\Omega)$, the Ritz projection $\mathcal R_hv\in V_h$ defined in~\eqref{eq:nh-ritz-def} satisfies
	\begin{align}
		\|\nabla(v-\mathcal R_hv)\|_{L^2(\Omega)}&\le C h^k|v|_{H^{k+1}(\Omega)}, \label{eq:Ritz-H-global}\\
		\|v-\mathcal R_hv\|_{L^2(\Gamma)}&\le C h^{k+1/2}|v|_{H^{k+1}(\Omega)}, \label{eq:Ritz-H-boundary}\\
		\|\partial_n(v-\mathcal R_hv)\|_{L^2(\Gamma)}&\le C h^{k-1/2}|v|_{H^{k+1}(\Omega)}.  \label{eq:Ritz-H-normal}
	\end{align}
\end{lemma}

\subsection{Energy estimate}

The next lemma records the part of the traditional analysis that works over the full penalty scale. For $\alpha\ge1$ it follows from the diagonal identity~\eqref{eq:diag-cancel}. For $\alpha<1$ it follows from \Cref{thm:weak-infsup}.

\begin{lemma}\label{lem:traditional-energy}
	Under Assumption~\ref{assump:standing}, let $\alpha\in\{-\infty\}\cup\mathbb R$; when $\alpha<1$, use meshes with $h<h_0$ as in \Cref{thm:weak-infsup}. Let $u\in H^{k+1}(\Omega)$ be the exact solution of~\eqref{org} and let $u_h\in V_h$ solve~\eqref{eq:prelim-nitsche}. Then
	\begin{align}
		\|\nabla(u-u_h)\|_{L^2(\Omega)}&\le C h^k|u|_{H^{k+1}(\Omega)}, \label{eq:traditional-energy-grad}\\
		\|u-u_h\|_{L^2(\Gamma)}&\le C h^{k+1/2}|u|_{H^{k+1}(\Omega)}.  \label{eq:traditional-energy-boundary}
	\end{align}
	Moreover, with $\alpha_\ast:=\max\{1,\alpha\}$,
	\begin{equation}\label{eq:traditional-Ru-boundary}
		\|\mathcal R_hu-u_h\|_{L^2(\Gamma)} \le C h^{k+\alpha_\ast/2}|u|_{H^{k+1}(\Omega)}. 
	\end{equation}
\end{lemma}

\begin{proof}
	Galerkin orthogonality gives
	\begin{equation}\label{eq:energy-residual-eq}
		a(\mathcal R_hu-u_h,v_h)=-a(u-\mathcal R_hu,v_h)\qquad\forall v_h\in V_h. 
	\end{equation}
	
	\smallskip\noindent\emph{Bound on $a(u-\mathcal R_h u,v_h)$.}  Expanding $a$ according to~\eqref{eq:ah-def},
	\begin{equation*}
		\begin{aligned}
			a(u-\mathcal R_hu,v_h)
			&= (\nabla(u-\mathcal R_hu),\nabla v_h)_\Omega
			- \langle\partial_n(u-\mathcal R_hu),v_h\rangle_\Gamma \\
			&\quad + \langle u-\mathcal R_hu,\partial_n v_h\rangle_\Gamma
			+ h^{-\alpha}\langle u-\mathcal R_hu,v_h\rangle_\Gamma. 
		\end{aligned}
	\end{equation*}
	Since $v_h|_\Gamma\in M_h$, the Ritz boundary orthogonality~\eqref{eq:Ritz-boundary-orthogonality} kills the penalty term for every $\alpha$ (including $\alpha=-\infty$, where the term is absent by convention).  The three remaining terms are bounded by Cauchy--Schwarz together with the Ritz estimates~\eqref{eq:Ritz-H-global}--\eqref{eq:Ritz-H-normal}, the norm definition~\eqref{eq:reference-nitsche-norm}, and the inverse trace~\eqref{eq:grad-inverse-trace}:
	\begin{align*}
		\bigl|(\nabla(u-\mathcal R_hu),\nabla v_h)_\Omega\bigr|
		&\le \|\nabla(u-\mathcal R_hu)\|_{L^2(\Omega)}\|\nabla v_h\|_{L^2(\Omega)} \le Ch^k|u|_{H^{k+1}(\Omega)}\vertiii{v_h}_h, \\
		\bigl|\langle\partial_n(u-\mathcal R_hu),v_h\rangle_\Gamma\bigr|
		&\le \|\partial_n(u-\mathcal R_hu)\|_{L^2(\Gamma)}\|v_h\|_{L^2(\Gamma)}\\
		&\le Ch^{k-1/2}|u|_{H^{k+1}(\Omega)}\cdot h^{\alpha_\ast/2}\vertiii{v_h}_h
		\le Ch^k|u|_{H^{k+1}(\Omega)}\vertiii{v_h}_h, \\
		\bigl|\langle u-\mathcal R_hu,\partial_n v_h\rangle_\Gamma\bigr|
		&\le \|u-\mathcal R_hu\|_{L^2(\Gamma)}\|\partial_n v_h\|_{L^2(\Gamma)}\\
		&\le Ch^{k+1/2}|u|_{H^{k+1}(\Omega)}\cdot Ch^{-1/2}\vertiii{v_h}_h
		= Ch^k|u|_{H^{k+1}(\Omega)}\vertiii{v_h}_h. 
	\end{align*}
	In the second estimate we used $\|v_h\|_{L^2(\Gamma)}\le h^{\alpha_\ast/2}\vertiii{v_h}_h$ from~\eqref{eq:reference-nitsche-norm} together with $\alpha_\ast\ge1$ and $h\le1$ to absorb the residual $h$-power; in the third, the inverse trace~\eqref{eq:grad-inverse-trace} summed over boundary elements gives $\|\partial_n v_h\|_{L^2(\Gamma)}\le Ch^{-1/2}\|\nabla v_h\|_{L^2(\Omega)}\le Ch^{-1/2}\vertiii{v_h}_h$.  Combining,
	\begin{equation}\label{eq:traditional-residual-fullscale}
		|a(u-\mathcal R_hu,v_h)|\le Ch^k|u|_{H^{k+1}(\Omega)}\vertiii{v_h}_h. 
	\end{equation}
	
	\smallskip\noindent\emph{Energy bound on $\mathcal R_h u-u_h$.}  If $\alpha\ge1$, take $v_h=\mathcal R_hu-u_h$ in~\eqref{eq:energy-residual-eq} and use the diagonal identity~\eqref{eq:diag-cancel} to obtain
	\begin{equation*}
		\vertiii{\mathcal R_hu-u_h}_h^2 = a(\mathcal R_hu-u_h,\mathcal R_hu-u_h) = -a(u-\mathcal R_hu,\mathcal R_hu-u_h). 
	\end{equation*}
	If $\alpha<1$, the inf--sup condition~\eqref{eq:weak-infsup} applied to $\mathcal R_hu-u_h$ together with~\eqref{eq:energy-residual-eq} gives the same conclusion up to the inf--sup constant $c$.  In either case,~\eqref{eq:traditional-residual-fullscale} yields
	\begin{equation}\label{eq:traditional-energy-Ru}
		\vertiii{\mathcal R_hu-u_h}_h \le Ch^k|u|_{H^{k+1}(\Omega)}. 
	\end{equation}
	
	\smallskip\noindent\emph{Proof of \eqref{eq:traditional-Ru-boundary}.}  From~\eqref{eq:reference-nitsche-norm}, $\|w\|_{L^2(\Gamma)}\le h^{\alpha_\ast/2}\vertiii{w}_h$.  Applying this with $w=\mathcal R_hu-u_h$ and using~\eqref{eq:traditional-energy-Ru} gives~\eqref{eq:traditional-Ru-boundary}.
	
	\smallskip\noindent\emph{Proof of \eqref{eq:traditional-energy-grad}.}  By the triangle inequality, the Ritz estimate~\eqref{eq:Ritz-H-global}, and~\eqref{eq:traditional-energy-Ru},
	\begin{equation*}
		\|\nabla(u-u_h)\|_{L^2(\Omega)}
		\le \|\nabla(u-\mathcal R_hu)\|_{L^2(\Omega)} + \vertiii{\mathcal R_hu-u_h}_h
		\le Ch^k|u|_{H^{k+1}(\Omega)}. 
	\end{equation*}
	
	\smallskip\noindent\emph{Proof of \eqref{eq:traditional-energy-boundary}.}  By the triangle inequality, the Ritz estimate~\eqref{eq:Ritz-H-boundary}, and~\eqref{eq:traditional-Ru-boundary},
	\begin{equation*}
		\|u-u_h\|_{L^2(\Gamma)}
		\le Ch^{k+1/2}|u|_{H^{k+1}(\Omega)} + Ch^{k+\alpha_\ast/2}|u|_{H^{k+1}(\Omega)}
		\le Ch^{k+1/2}|u|_{H^{k+1}(\Omega)}.
	\end{equation*}
\end{proof}

\subsection{Duality}\label{sec:traditional-duality}

The energy estimate is now combined with a standard Aubin--Nitsche duality argument to control $\|u-u_h\|_{L^2(\Omega)}$ on the full penalty scale.

\begin{lemma}\label{lem:standard-duality-T3}
	Under Assumption~\ref{assump:standing}, let $\alpha\in\{-\infty\}\cup\mathbb R$; when $\alpha<1$, use meshes with $h<h_0$ as in \Cref{thm:weak-infsup}. Let $u\in H^{k+1}(\Omega)$ be the exact solution of~\eqref{org}, let $u_h\in V_h$ solve~\eqref{eq:prelim-nitsche}, and set $\alpha_\ast:=\max\{1,\alpha\}$. Then
	\begin{equation}\label{eq:traditional-final-rate}
		\|u-u_h\|_{L^2(\Omega)}
		\le C\bigl(h^{k+1}+h^{k+\alpha_\ast/2}\bigr)|u|_{H^{k+1}(\Omega)}. 
	\end{equation}
\end{lemma}

\begin{proof}
	Let $\Psi\in H^2(\Omega)\cap H_0^1(\Omega)$ solve the dual problem $-\Delta\Psi=u-u_h$ in $\Omega$ with $\Psi=0$ on $\Gamma$; by convex-domain elliptic regularity~\eqref{dual},
	\begin{equation}\label{eq:standard-dual-regularity}
		\|\Psi\|_{H^2(\Omega)}\le C\|u-u_h\|_{L^2(\Omega)}. 
	\end{equation}
	
	\emph{Dual identity.}  Since $I_h\Psi\in V_h^0$, integration by parts together with Galerkin orthogonality~\eqref{eq:galerkin-orthogonality} give
	\[
	\|u-u_h\|_{L^2(\Omega)}^2
	= a(u-u_h,\Psi-I_h\Psi) - 2\langle u-u_h,\partial_n\Psi\rangle_\Gamma,
	\]
	Decompose $\partial_n\Psi=\Pi_h^\Gamma\partial_n\Psi+(I-\Pi_h^\Gamma)\partial_n\Psi$.  Since $\Pi_h^\Gamma\partial_n\Psi\in M_h$, the boundary orthogonality~\eqref{eq:Ritz-boundary-orthogonality} of the projected-boundary Ritz projection $\mathcal R_h$ converts the projected part to its Ritz form, $\langle u-u_h,\Pi_h^\Gamma\partial_n\Psi\rangle_\Gamma = \langle \mathcal R_hu-u_h,\Pi_h^\Gamma\partial_n\Psi\rangle_\Gamma$, and the dual identity reads
	\begin{equation}\label{eq:standard-duality-T123}
		\|u-u_h\|_{L^2(\Omega)}^2 = T_1-2T_2-2T_3,
	\end{equation}
	where
	\begin{align*}
		T_1&:=a(u-u_h,\Psi-I_h\Psi), \\
		T_2&:=\langle u-u_h,(I-\Pi_h^\Gamma)\partial_n\Psi\rangle_\Gamma,
	\end{align*}
	and the projected boundary term, taken in its Ritz form, is
	\begin{equation}\label{eq:T3-Ritz-form}
		T_3 := \langle \mathcal R_hu-u_h,\Pi_h^\Gamma\partial_n\Psi\rangle_\Gamma. 
	\end{equation}

	\emph{Universal $T_1+T_2$ bound.}  Since $(\Psi-I_h\Psi)|_\Gamma=0$, the penalty term and the term containing $\partial_n(u-u_h)$ in $T_1$ vanish, leaving
	\[
	T_1=(\nabla(u-u_h),\nabla(\Psi-I_h\Psi))_\Omega+\langle u-u_h,\partial_n(\Psi-I_h\Psi)\rangle_\Gamma. 
	\]
	\Cref{lem:traditional-energy}, the interpolation estimate~\eqref{eq:H-interp-volume}, and the normal-trace interpolation from~\eqref{eq:scaled-trace} give
	\[
	|T_1|\le Ch^{k+1}|u|_{H^{k+1}(\Omega)}\|\Psi\|_{H^2(\Omega)}. 
	\]
	For $T_2$, the projection estimate~\eqref{eq:normal-Pi-est} and the boundary error bound~\eqref{eq:traditional-energy-boundary} give the same bound $|T_2|\le Ch^{k+1}|u|_{H^{k+1}(\Omega)}\|\Psi\|_{H^2(\Omega)}$, so
	\begin{equation}\label{eq:T1-T2-universal}
		|T_1|+|T_2| \le Ch^{k+1}|u|_{H^{k+1}(\Omega)}\|\Psi\|_{H^2(\Omega)}. 
	\end{equation}
	This bound is independent of $\alpha$ and is the only ingredient of the duality argument reused by \Cref{thm:main-L2,thm:weak-L2}.
	
	\emph{Traditional $T_3$ bound.}  Using the Ritz form~\eqref{eq:T3-Ritz-form}, Cauchy--Schwarz, the $L^2(\Gamma)$ stability of $\Pi_h^\Gamma$, the trace theorem~\eqref{trace_theorem}, and~\eqref{eq:traditional-Ru-boundary} give
	\begin{equation}\label{eq:T3-projected-traditional}
		|T_3|\le \|\mathcal R_hu-u_h\|_{L^2(\Gamma)}\|\Pi_h^\Gamma\partial_n\Psi\|_{L^2(\Gamma)}
		\le C h^{k+\alpha_\ast/2}|u|_{H^{k+1}(\Omega)}\|\Psi\|_{H^2(\Omega)}. 
	\end{equation}
	
	\emph{Conclusion.}  Substituting~\eqref{eq:T1-T2-universal} and~\eqref{eq:T3-projected-traditional} into~\eqref{eq:standard-duality-T123}, using~\eqref{eq:standard-dual-regularity}, and dividing by $\|u-u_h\|_{L^2(\Omega)}$ proves~\eqref{eq:traditional-final-rate}.
\end{proof}

\begin{remark}\label{rem:half-order-loss}
	At the classical stabilized scaling $\alpha=1$ and at Burman's penalty-free endpoint $\alpha=-\infty$, we have $\alpha_\ast=1$, so~\eqref{eq:traditional-final-rate} reduces to
	\[
	\|u-u_h\|_{L^2(\Omega)}\le C h^{k+1/2}|u|_{H^{k+1}(\Omega)}. 
	\]
	Both formulations therefore incur the same half-order loss compared with the optimal rate $h^{k+1}$ that holds for adjoint-consistent methods.
\end{remark}

\subsection{Discussion}\label{sec:obstruction}

Remark~\ref{rem:half-order-loss} predicts the same half-order rate $h^{k+1/2}$ for both the classical stabilized method ($\alpha=1$) and Burman's penalty-free method ($\alpha=-\infty$). Numerical experiments on smooth test problems, however, consistently produce the full $h^{k+1}$ rate; the tables in Section~\ref{sec:numerics} reproduce this behavior. This discrepancy between theory and computation has long persisted: attempts to construct an example that realizes the predicted half-order loss had repeatedly failed (cf.\ the remark of Burman quoted in Section~\ref{sec:introduction}).  This created a natural ambiguity: is the half-order gap a defect of the nonsymmetric Nitsche method, or only a defect of the classical analysis?

In this paper we resolve this ambiguity. Section~\ref{sec:numerics} constructs a parametric family of manufactured solutions at the rough endpoint of $H^{k+1}$ regularity for which the classical method ($\alpha=1$) and Burman's penalty-free method ($\alpha=-\infty$) both clearly exhibit the predicted half-order loss, attaining the rate $h^{k+1/2}$. The estimate of \Cref{lem:standard-duality-T3} is therefore sharp under merely $H^{k+1}$ regularity, and there is no room to improve the theory at that level. This does not, however, explain why the same penalty scalings produce the optimal rate $h^{k+1}$ in smooth tests. A sufficiently strong penalty $\alpha\ge 2$ restores optimality in the estimate~\eqref{eq:traditional-final-rate}, but this elementary repair does not address the central phenomenon: the optimal rates observed in computations at $\alpha=1$ and $\alpha=-\infty$.

The phenomenon of optimal convergence in smooth tests motivates measuring regularity at the endpoint $W^{k+1,\infty}(\Omega)$. To unify the analysis along the full Sobolev scale, we work with $W^{k+1,p}(\Omega)$ for $p\in[2,+\infty]$. Sections~\ref{sec:main} and~\ref{sec:weak-penalty} prove the unified estimate of \Cref{thm:unified-L2}:
\begin{equation}\label{eq:story-main-rate-preview}
	\begin{split}
		\|u-u_h\|_{L^2(\Omega)} &\le C h^{r}|u|_{W^{k+1,p}(\Omega)},\quad  r = \min\bigl\{k+1,\,k+\max\{1,\alpha\}-1/p\bigr\}.
	\end{split}
\end{equation}
Specializing~\eqref{eq:story-main-rate-preview} to the classical stabilized method ($\alpha=1$) and Burman's penalty-free method ($\alpha=-\infty$) gives the same rate for both, since $\max\{1,\alpha\}=1$ in either case; at the two endpoints of the Sobolev scale,
\[
\|u-u_h\|_{L^2(\Omega)}\le
\begin{cases}
	C h^{k+1/2}\,|u|_{H^{k+1}(\Omega)}, & u\in H^{k+1}(\Omega),\\[2pt]
	C h^{k+1}\,|u|_{W^{k+1,\infty}(\Omega)}, & u\in W^{k+1,\infty}(\Omega).
\end{cases}
\]

These two endpoint rates together resolve the long-standing ambiguity: the half-order loss is genuine under merely $H^{k+1}$ regularity, in agreement with Remark~\ref{rem:half-order-loss}, while the optimal rate is recovered under $W^{k+1,\infty}$ regularity, explaining the convergence observed in smooth tests.

\section{Sharp \texorpdfstring{$L^2$}{L2} estimate for \texorpdfstring{$\alpha\ge 1$}{α≥1}}\label{sec:main}

As motivated by the discussion in Section~\ref{sec:obstruction}, we now sharpen the analysis under the finer regularity scale $u\in W^{k+1,p}(\Omega)$, $p\in[2,+\infty]$.  \Cref{lem:standard-duality-T3} provides the universal bound~\eqref{eq:T1-T2-universal}, which controls $T_1$ and $T_2$ at the optimal order $h^{k+1}$ even at the rough endpoint $u\in H^{k+1}(\Omega)$, so the only term to be improved is the projected boundary component, taken in its Ritz form~\eqref{eq:T3-Ritz-form},
\[
T_3=\langle \mathcal R_hu-u_h,\Pi_h^\Gamma\partial_n\Psi\rangle_\Gamma,
\]
where $\Psi\in H^2(\Omega)\cap H_0^1(\Omega)$ solves the dual problem $-\Delta\Psi=u-u_h$ in $\Omega$, $\Psi=0$ on $\Gamma$.  The remainder of this section focuses on a sharper estimate of $T_3$ for $\alpha\ge 1$.

A direct Cauchy--Schwarz argument fails to yield the optimal rate, as we now show. Even with the stronger regularity $u\in W^{k+1,p}(\Omega)$, a direct estimate via Cauchy--Schwarz, the boundary control of the energy estimate, and H\"older's inequality $|u|_{H^{k+1}(\Omega)}\le|\Omega|^{1/2-1/p}|u|_{W^{k+1,p}(\Omega)}$ gives
\begin{align}\label{eq:direct-Wkp-still-loss}
	\begin{split}
		|T_3|
		&\le \|\mathcal R_hu-u_h\|_{L^2(\Gamma)}\|\Pi_h^\Gamma\partial_n\Psi\|_{L^2(\Gamma)}
		\le C h^{k+\alpha/2}|u|_{H^{k+1}(\Omega)}\|\Psi\|_{H^2(\Omega)} \\
		&\le C|\Omega|^{1/2-1/p}\,h^{k+\alpha/2}|u|_{W^{k+1,p}(\Omega)}\|\Psi\|_{H^2(\Omega)}.
	\end{split}
\end{align}
The $W^{k+1,p}$ regularity contributes only the mesh-independent constant $|\Omega|^{1/2-1/p}$; for $\alpha=1$ the rate is still $h^{k+1/2}$, independently of $p$.

The obstruction is structural: Cauchy--Schwarz can only extract the square root $h^{\alpha/2}$ from the boundary energy norm. The full exponent $h^\alpha$, however, already lives inside the bilinear form $a$ as the multiplier of its penalty term $h^{-\alpha}\langle u-u_h,v_h|_\Gamma\rangle_\Gamma$, which has the same boundary-inner-product structure as $T_3$. Choosing $v_h\in V_h$ with $v_h|_\Gamma=\Pi_h^\Gamma\partial_n\Psi$ (admissible since $\Pi_h^\Gamma\partial_n\Psi\in M_h$) makes the penalty contribution exactly $h^{-\alpha}T_3$ by boundary orthogonality of $\mathcal R_h$, as we now verify.  The boundary orthogonality~\eqref{eq:Ritz-boundary-orthogonality} of $\mathcal R_h$ with the test function $\mu_h=\Pi_h^\Gamma\partial_n\Psi\in M_h$ reads $\langle u-\mathcal R_hu,\Pi_h^\Gamma\partial_n\Psi\rangle_\Gamma=0$, so the Ritz form~\eqref{eq:T3-Ritz-form} of $T_3$ extends to
\[
T_3 = \langle \mathcal R_hu-u_h,\Pi_h^\Gamma\partial_n\Psi\rangle_\Gamma
= \langle u-u_h,\Pi_h^\Gamma\partial_n\Psi\rangle_\Gamma. 
\]
With $v_h|_\Gamma=\Pi_h^\Gamma\partial_n\Psi$, this identifies $\langle u-u_h,v_h\rangle_\Gamma=T_3$, so the penalty term $h^{-\alpha}\langle u-u_h,v_h\rangle_\Gamma$ equals $h^{-\alpha}T_3$ exactly.  Galerkin orthogonality $a(u-u_h,v_h)=0$ then yields the exact identity, which we refer to as the \emph{projected boundary functional identity}:
\begin{align}\label{eq:boundary-functional-identity}
	h^{-\alpha}T_3
	&=-(\nabla(u-u_h),\nabla v_h)_\Omega
	+\langle\partial_n(u-u_h),v_h\rangle_\Gamma  -\langle u-u_h,\partial_n v_h\rangle_\Gamma. 
\end{align}
This identity converts the Cauchy--Schwarz factor $h^{\alpha/2}$ on the boundary term $T_3$ into the full penalty factor $h^\alpha$ multiplying volume and trace residuals. Specializing \Cref{cor:dual-lift-scales} to $z=\Psi$, the lifting $v_h$ has the boundary-layer scaling
\begin{equation}\label{eq:opening-lift-scales}
	\|\nabla v_h\|_{L^2(\Omega)}\le Ch^{-1/2}\|\Psi\|_{H^2(\Omega)},
	\qquad
	\|\partial_n v_h\|_{L^2(\Gamma)}\le Ch^{-1}\|\Psi\|_{H^2(\Omega)}. 
\end{equation}

Estimating the right-hand side of~\eqref{eq:boundary-functional-identity} using the global energy estimate of \Cref{lem:traditional-energy} and the lifting scalings~\eqref{eq:opening-lift-scales}, each of the three residual terms is bounded by $Ch^{k-1/2}|u|_{H^{k+1}(\Omega)}\|\Psi\|_{H^2(\Omega)}$; multiplying through by $h^\alpha$ gives
\[
|T_3|\le Ch^{k+\alpha-1/2}|u|_{H^{k+1}(\Omega)}\|\Psi\|_{H^2(\Omega)}. 
\]
This represents a genuine improvement over~\eqref{eq:direct-Wkp-still-loss}: the rate sharpens from $h^{k+\alpha/2}$ to $h^{k+\alpha-1/2}$. A factor $h^{1/2-1/p}$ is still missing; replacing $H^{k+1}$ by $W^{k+1,p}$ globally contributes only the mesh-independent constant $|\Omega|^{1/2-1/p}$ in~\eqref{eq:direct-Wkp-still-loss}.

To extract such a factor from H\"older's inequality, we must apply it on a set of measure $O(h)$ rather than on the whole domain. Boundary-strip localization, a well-established technique with a long history~\cite{MR1974490,MR3319574,PfeffererWinkler2019,MR4057429}, provides exactly such a set. We therefore choose a special lifting $v_h$ that vanishes outside the boundary strip $\Omega_h^\Gamma$ defined in~\eqref{eq:strip-def}.
With this choice, the volume residual $(\nabla(u-u_h),\nabla v_h)_\Omega$ in~\eqref{eq:boundary-functional-identity} reduces to the strip integral $(\nabla(u-u_h),\nabla v_h)_{\Omega_h^\Gamma}$, and H\"older's inequality on $\Omega_h^\Gamma$ now supplies precisely the missing factor $h^{1/2-1/p}$.

Accordingly, the three quantities to be sharpened are
\begin{align}\label{eq:critical-strip-error-opening}
	\begin{split}
		&\|\nabla(u-u_h)\|_{L^2(\Omega_h^\Gamma)}
		+h^{1/2}\|\partial_n(u-u_h)\|_{L^2(\Gamma)}
		+h^{-1/2}\|u-u_h\|_{L^2(\Gamma)} \\
		&\qquad\le C h^{k+1/2-1/p}|u|_{W^{k+1,p}(\Omega)}. 
	\end{split}
\end{align}
Multiplying through by $h^\alpha$ then yields the sharp bound
\[
|T_3|\le Ch^{k+\alpha-1/p}|u|_{W^{k+1,p}(\Omega)}\|\Psi\|_{H^2(\Omega)},
\]
in which the favorable factor $h^{k+\alpha-1/p}$ replaces the half-order $h^{k+\alpha/2}$ of the direct approach. Together with the optimal bounds on $T_1$ and $T_2$ from \Cref{lem:standard-duality-T3} and the dual regularity $\|\Psi\|_{H^2}\le C\|u-u_h\|_{L^2}$, this delivers the sharp $L^2$ estimate
\[
\|u-u_h\|_{L^2(\Omega)}\le C h^{\min\{k+1,\,k+\alpha-1/p\}}|u|_{W^{k+1,p}(\Omega)},
\]
the main result of this section.

The proof of~\eqref{eq:critical-strip-error-opening} is the technical core of this section and the proof proceeds in four steps: Step~1 establishes the $W^{k+1,p}$ interpolation and projection estimates on the strip and on $\Gamma$; Step~2 constructs the boundary-supported lifting and proves the boundary-layer Ritz estimate; Step~3 transfers these estimates to the error near the boundary; Step~4 returns to the dual identity~\eqref{eq:boundary-functional-identity} and uses it to control $T_3$.

\subsection{\texorpdfstring{$W^{k+1,p}$}{Wk+1,p} interpolation and projection estimates}

We collect the $W^{k+1,p}$ approximation estimates for the Lagrange interpolant $I_h$ (\Cref{lem:delta-Wkp}) and the boundary $L^2$ projection $\Pi_h^\Gamma$ (\Cref{lem:boundary-L2-proj}).

\begin{lemma}\label{lem:delta-Wkp}
	Under Assumption~\ref{assump:standing}, let $v\in W^{k+1,p}(\Omega)$ for some $p\in[2,\infty]$.  Then
	\begin{align}
		\|v-I_hv\|_{L^2(\Gamma)} &\le C h^{k+1-1/p}|v|_{W^{k+1,p}(\Omega)}, \label{eq:lem-delta-Wkp}\\
		\|\nabla(v-I_hv)\|_{L^2(\Gamma)} &\le C h^{k-1/p}|v|_{W^{k+1,p}(\Omega)}, \label{eq:interp-tang}\\
		\|\nabla(v-I_hv)\|_{L^2(\Omega_h^\Gamma)} &\le C h^{k+1/2-1/p}|v|_{W^{k+1,p}(\Omega)}.  \label{eq:interp-strip-grad}
	\end{align}
\end{lemma}

\begin{proof}
	The standard $L^p$ Lagrange interpolation estimate~\cite[Theorem~3.1.6]{Ciarlet1978} combined with H\"older's inequality on each element gives the local bound
	\begin{equation}\label{eq:cross-scale-interp}
		\|D^j(v-I_hv)\|_{L^2(T)}\le C h_T^{d(1/2-1/p)+k+1-j}|v|_{W^{k+1,p}(T)},\qquad 0\le j\le k+1,\ T\in\mathcal T_h,
	\end{equation}
	with the convention $1/\infty=0$ when $p=\infty$.
	
	\smallskip\noindent\emph{Boundary estimates~\eqref{eq:lem-delta-Wkp}--\eqref{eq:interp-tang}.}  For each boundary face $E\in\mathcal E_h^\Gamma$ with parent element $T_E$, the scaled trace inequality~\eqref{eq:scaled-trace} applied componentwise to $D^m(v-I_hv)$ gives
	\begin{align*}
		\|D^m(v-I_hv)\|_{L^2(E)}^2
		&\le C h_{T_E}^{-1}\|D^m(v-I_hv)\|_{L^2(T_E)}^2 \\
		&\quad + C h_{T_E}\|D^{m+1}(v-I_hv)\|_{L^2(T_E)}^2,\quad m=0,1. 
	\end{align*}
	Substituting~\eqref{eq:cross-scale-interp} (at $j=m$ and $j=m+1$) and using $h_{T_E}\simeq h$,
	\begin{equation*}
		\|D^m(v-I_hv)\|_{L^2(E)}^2 \le C h^{2(k+1-m)+2d(1/2-1/p)-1}|v|_{W^{k+1,p}(T_E)}^2. 
	\end{equation*}
	Summing over $E\in\mathcal E_h^\Gamma$ (each $T\in\mathcal T_h^\Gamma$ has a uniformly bounded number of boundary faces) yields
	\begin{equation*}
		\|D^m(v-I_hv)\|_{L^2(\Gamma)}^2
		\le C h^{2(k+1-m)+2d(1/2-1/p)-1}\sum_{T\in\mathcal T_h^\Gamma}|v|_{W^{k+1,p}(T)}^2. 
	\end{equation*}
	The discrete H\"older inequality on the $\#\mathcal T_h^\Gamma\sim h^{-(d-1)}$ boundary elements gives, uniformly in $p\in[2,\infty]$ with the convention $1/\infty=0$,
	\begin{align}
		\begin{split}
			\sum_{T\in\mathcal T_h^\Gamma}|v|_{W^{k+1,p}(T)}^2
			&\le \bigl(\#\mathcal T_h^\Gamma\bigr)^{1-2/p}\Bigl(\sum_{T\in\mathcal T_h^\Gamma}|v|_{W^{k+1,p}(T)}^p\Bigr)^{2/p} \\
			&\le C h^{-(d-1)(1-2/p)}|v|_{W^{k+1,p}(\Omega)}^2. 
		\end{split}
	\end{align}
	Combining,
	\begin{equation*}
		\|D^m(v-I_hv)\|_{L^2(\Gamma)}^2 \le C h^{2(k+1-m-1/p)}|v|_{W^{k+1,p}(\Omega)}^2,
	\end{equation*}
	which gives~\eqref{eq:lem-delta-Wkp} at $m=0$ and~\eqref{eq:interp-tang} at $m=1$.  The bound on the elementwise normal derivative follows from $|\partial_n w|=|(\nabla w)\cdot n|\le|\nabla w|$ pointwise on each $E$.
	
	\smallskip\noindent\emph{Strip estimate~\eqref{eq:interp-strip-grad}.}  H\"older's inequality on $\Omega_h^\Gamma$ together with the measure bound $|\Omega_h^\Gamma|\le Ch$ and the standard $L^p$ Lagrange estimate gives
	\begin{align*}
		\|\nabla(v-I_hv)\|_{L^2(\Omega_h^\Gamma)} &
		\le |\Omega_h^\Gamma|^{1/2-1/p}\|\nabla(v-I_hv)\|_{L^p(\Omega_h^\Gamma)}\\
		&\le Ch^{1/2-1/p}\cdot Ch^k|v|_{W^{k+1,p}(\Omega)}.
	\end{align*}
\end{proof}

The next lemma transfers these bounds to the boundary $L^2$ projection $\Pi_h^\Gamma$.
\begin{lemma}\label{lem:boundary-L2-proj}
	Under Assumption~\ref{assump:standing}, let $v\in W^{k+1,p}(\Omega)$ for some $p\in[2,\infty]$.  Then
	\begin{align}
		\|v|_\Gamma-\Pi_h^\Gamma(v|_\Gamma)\|_{L^2(\Gamma)}
		&\le C h^{k+1-1/p}|v|_{W^{k+1,p}(\Omega)}, \label{eq:Pi-Wkp-boundary}\\
		\|v|_\Gamma-\Pi_h^\Gamma(v|_\Gamma)\|_{H^{1/2}(\Gamma)}
		&\le C h^{k+1/2-1/p}|v|_{W^{k+1,p}(\Omega)}.  \label{eq:Pi-Hhalf-Wkp}
	\end{align}
\end{lemma}

\begin{proof}
	For~\eqref{eq:Pi-Wkp-boundary}: since $\Pi_h^\Gamma(v|_\Gamma)$ is the $L^2(\Gamma)$ best approximation of $v|_\Gamma$ in $M_h$ and $(I_hv)|_\Gamma\in M_h$,
	\[
	\|v|_\Gamma-\Pi_h^\Gamma(v|_\Gamma)\|_{L^2(\Gamma)}
	\le \|v-I_hv\|_{L^2(\Gamma)}
	\le Ch^{k+1-1/p}|v|_{W^{k+1,p}(\Omega)},
	\]
	by~\eqref{eq:lem-delta-Wkp}.
	For~\eqref{eq:Pi-Hhalf-Wkp}, we first establish the $H^1(\Gamma)$ bound
	\begin{equation}\label{eq:Pi-H1-Wkp}
		\|v|_\Gamma-\Pi_h^\Gamma(v|_\Gamma)\|_{H^1(\Gamma)} \le Ch^{k-1/p}|v|_{W^{k+1,p}(\Omega)},
	\end{equation}
	and then interpolate.  Decompose
	\[
	v|_\Gamma-\Pi_h^\Gamma(v|_\Gamma) = (v-I_hv)|_\Gamma + \bigl((I_hv)|_\Gamma-\Pi_h^\Gamma(v|_\Gamma)\bigr). 
	\]
	For the first term, the surface gradient is the tangential projection of the ambient gradient, so $|\nabla_\Gamma w|\le|\nabla w|$ pointwise on each face for the trace $w=(v-I_hv)|_\Gamma$ of the bulk function $v-I_hv$; together with~\eqref{eq:lem-delta-Wkp}--\eqref{eq:interp-tang},
	\[
	\|(v-I_hv)|_\Gamma\|_{H^1(\Gamma)} \le Ch^{k-1/p}|v|_{W^{k+1,p}(\Omega)}. 
	\]
	The second term lies in $M_h$, so applying the standard polynomial inverse estimate~\cite[Lemma~4.5.3, Eq.~(4.5.4), p.~111]{Brenner2008} on each boundary face $E\in\mathcal E_h^\Gamma$ (a shape-regular $(d-1)$-dimensional simplex) and summing using quasi-uniformity gives $\|\mu_h\|_{H^1(\Gamma)}\le Ch^{-1}\|\mu_h\|_{L^2(\Gamma)}$ for $\mu_h\in M_h$; combined with the triangle inequality, \eqref{eq:lem-delta-Wkp}, and~\eqref{eq:Pi-Wkp-boundary},
	\begin{align*}
		\|(I_hv)|_\Gamma-\Pi_h^\Gamma(v|_\Gamma)\|_{H^1(\Gamma)}
		&\le Ch^{-1}\bigl(\|v-I_hv\|_{L^2(\Gamma)}+\|v|_\Gamma-\Pi_h^\Gamma(v|_\Gamma)\|_{L^2(\Gamma)}\bigr) \\
		&\le Ch^{k-1/p}|v|_{W^{k+1,p}(\Omega)}. 
	\end{align*}
	Adding the two bounds gives~\eqref{eq:Pi-H1-Wkp}.  The Sobolev interpolation inequality~\cite[Chapter~1, Proposition~2.3]{Lions1972}
	\[
	\|w\|_{H^{1/2}(\Gamma)} \le C\|w\|_{L^2(\Gamma)}^{1/2}\|w\|_{H^1(\Gamma)}^{1/2},
	\]
	applied with $w=v|_\Gamma-\Pi_h^\Gamma(v|_\Gamma)$ combines~\eqref{eq:Pi-Wkp-boundary} and~\eqref{eq:Pi-H1-Wkp} into~\eqref{eq:Pi-Hhalf-Wkp}.
\end{proof}

\subsection{Boundary-layer Ritz estimate}

The proof of the boundary-layer Ritz estimate (\Cref{lem:nh-ritz-boundary-layer} below) uses the boundary-supported lifting $E_h$ of \Cref{lem:disc-lift} together with the pointwise gradient stability of the homogeneous Ritz projection $\mathcal R_h^0:H_0^1(\Omega)\to V_h^0$, defined by
\begin{equation}
	(\nabla\mathcal R_h^0w,\nabla w_h)_\Omega=(\nabla w,\nabla w_h)_\Omega\qquad\forall w_h\in V_h^0.
\end{equation}
The pointwise gradient stability estimate
\begin{equation}\label{eq:DRS-pointwise-bound}
	|\nabla\mathcal R_h^0w(x)|\le C\mathcal M(|\nabla w|)(x)\qquad\text{for a.e. }x\in\Omega,
\end{equation}
is due to Diening--Rolfes--Salgado~\cite[Theorem~3.1]{DieningRolfesSalgado2024}, which applies to convex polytopes under shape-regular quasi-uniform meshes (Assumption~\ref{assump:standing}); here $\mathcal M$ is the Hardy--Littlewood maximal operator. Whenever \(\mathcal M\) is applied to a function defined on \(\Omega\),
we extend that function by zero to \(\mathbb R^d\).
We shall use the boundedness of \(\mathcal M\) on \(L^p(\mathbb R^d)\)
for \(2\le p<\infty\), and the elementary endpoint bound
\[
\|\mathcal M f\|_{L^\infty(\mathbb R^d)}
\le
\|f\|_{L^\infty(\mathbb R^d)} .
\]

\begin{lemma}\label{lem:nh-ritz-boundary-layer}
	Under Assumption~\ref{assump:standing}, let $v\in W^{k+1,p}(\Omega)$ for some $p\in[2,\infty]$, and let $\mathcal R_hv$ be defined by~\eqref{eq:nh-ritz-def}.  Then
	\begin{align}
		\|\nabla(v-\mathcal R_hv)\|_{L^2(\Omega_h^\Gamma)}
		&+h^{-1/2}\|v-\mathcal R_hv\|_{L^2(\Gamma)}
		+h^{1/2}\|\partial_n(v-\mathcal R_hv)\|_{L^2(\Gamma)} \notag\\
		&\le C h^{k+1/2-1/p}|v|_{W^{k+1,p}(\Omega)}.  \label{eq:Ritz-boundary-layer}
	\end{align}
\end{lemma}
\begin{proof}
	{
		We decompose the error $\rho := v - \mathcal R_hv$ into manageable components and estimate each part in the boundary layer $\Omega_h^\Gamma$ and on the boundary $\Gamma$.
		
		Let $\eta_0 := v - I_hv$ be the standard interpolation error. We introduce a modified interpolant $J_hv \in V_h$:
		\[
		J_hv := I_hv + E_h\bigl(\Pi_h^\Gamma(\eta_0|_\Gamma)\bigr).
		\]
		Using the trace property $(E_h\mu_h)|_\Gamma=\mu_h$ from~\eqref{eq:disc-lift-trace} and the fact that $\Pi_h^\Gamma$ preserves $(I_hv)|_\Gamma$, we verify that $J_hv|_\Gamma = \Pi_h^\Gamma(v|_\Gamma) = (\mathcal R_hv)|_\Gamma$. We then define:
		\[
		\psi_E := J_hv - I_hv, \quad \eta := v - J_hv = \eta_0 - \psi_E, \quad z_h := \mathcal R_hv - J_hv \in V_h^0.
		\]
		The total error is thus $\rho = v - \mathcal R_hv = \eta_0 - \psi_E - z_h$. By the $L^2(\Gamma)$ stability of $\Pi_h^\Gamma$ and the properties of the discrete lift $E_h$, we obtain the bounds for $\psi_E$:
		\begin{align}\label{eq:psiE-bounds}
			\begin{split}
				\|\nabla\psi_E\|_{L^2(\Omega)} \le Ch^{k+1/2-1/p}|v|_{W^{k+1,p}(\Omega)}, \\
				\|\partial_n\psi_E\|_{L^2(\Gamma)} \le Ch^{k-1/p}|v|_{W^{k+1,p}(\Omega)}.
			\end{split}
		\end{align}
	}

	To apply the pointwise estimate~\eqref{eq:DRS-pointwise-bound}, we identify $z_h$ with the homogeneous Ritz projection of a function in $H_0^1(\Omega)$. Let $q\in H^1(\Omega)$ be the harmonic extension of $\eta|_\Gamma$:
	\[
	q|_\Gamma = \eta|_\Gamma,\quad (\nabla q,\nabla\varphi)_\Omega = 0 \quad \forall\varphi\in H_0^1(\Omega). 
	\]
	Since $\eta|_\Gamma = v|_\Gamma - \Pi_h^\Gamma(v|_\Gamma)$, the trace theorem and~\eqref{eq:Pi-Hhalf-Wkp} give
	\begin{equation}\label{eq:q-energy-bound}
		\|\nabla q\|_{L^2(\Omega)}\le C\|\eta|_\Gamma\|_{H^{1/2}(\Gamma)}\le Ch^{k+1/2-1/p}|v|_{W^{k+1,p}(\Omega)}. 
	\end{equation}
	{
		Define $z := \eta - q \in H_0^1(\Omega)$. For any $w_h \in V_h^0$, the harmonicity of $q$ and the Ritz orthogonality~\eqref{eq:nh-ritz-def} imply $(\nabla z_h, \nabla w_h)_\Omega = (\nabla \eta, \nabla w_h)_\Omega = (\nabla z, \nabla w_h)_\Omega$. Thus, $z_h = \mathcal R_h^0 z$ is the homogeneous Ritz projection of $z$.
	}

	Applying the pointwise bound~\eqref{eq:DRS-pointwise-bound} to $z_h = \mathcal R_h^0 z$ and using the sublinearity of the Hardy-Littlewood maximal operator $\mathcal M$, we have
	\[
	|\nabla z_h| \le C\bigl(\mathcal M(|\nabla\eta_0|) + \mathcal M(|\nabla\psi_E|) + \mathcal M(|\nabla q|)\bigr) \quad \text{in } \Omega_h^\Gamma.
	\]
	For the first term, Hölder's inequality on \(\Omega_h^\Gamma\)
	(with measure \(|\Omega_h^\Gamma|\le Ch\)), the \(L^p(\mathbb R^d)\)
	boundedness of \(\mathcal M\) for \(2\le p<\infty\), the endpoint bound
	\(\|\mathcal M f\|_{L^\infty(\mathbb R^d)}
	\le \|f\|_{L^\infty(\mathbb R^d)}\) when \(p=\infty\), and standard
	interpolation estimates give
	\begin{equation*}
		\begin{split}
			\|\mathcal{M}(|\nabla\eta_0|)\|_{L^2(\Omega_\Gamma^h)}&\le Ch^{1/2-1/p}\|\mathcal{M}(|\nabla\eta_0|)\|_{L^p(\Omega_\Gamma^h)}\le Ch^{1/2-1/p}\|\mathcal{M} (|\nabla\eta_0|)\|_{L^p(\mathbb R^d)}\\
			&\le Ch^{1/2-1/p}\|\nabla\eta_0\|_{L^p(\Omega)}\le Ch^{k+1/2-1/p}|v|_{W^{k+1,p}(\Omega)}.
		\end{split}
	\end{equation*}
	For the second and third terms, the $L^2(\mathbb R^d)$-boundedness of $\mathcal M$ combined with~\eqref{eq:psiE-bounds} and~\eqref{eq:q-energy-bound} yields the identical bound. Consequently,
	\begin{equation}\label{eq:zh-bdry-layer}
		\|\nabla z_h\|_{L^2(\Omega_h^\Gamma)} \le Ch^{k+1/2-1/p}|v|_{W^{k+1,p}(\Omega)}. 
	\end{equation}
	A direct application of the strip estimate~\eqref{eq:interp-strip-grad} (or the same H\"older--Lagrange chain without $\mathcal M$) gives the same bound for $\|\nabla\eta_0\|_{L^2(\Omega_h^\Gamma)}$. Applying the triangle inequality to $\nabla\rho = \nabla\eta_0 - \nabla\psi_E - \nabla z_h$ and using~\eqref{eq:psiE-bounds} and~\eqref{eq:zh-bdry-layer}, we conclude
	\begin{equation}\label{eq:rho-bdry-layer-grad}
		\|\nabla\rho\|_{L^2(\Omega_h^\Gamma)} \le Ch^{k+1/2-1/p}|v|_{W^{k+1,p}(\Omega)}. 
	\end{equation}

	Since $z_h|_\Gamma = 0$, the boundary trace of the error is simply $\rho|_\Gamma = v|_\Gamma - \Pi_h^\Gamma(v|_\Gamma)$. The approximation property~\eqref{eq:Pi-Wkp-boundary} directly yields
	\begin{equation}\label{eq:rho-trace}
		h^{-1/2}\|\rho\|_{L^2(\Gamma)} \le Ch^{k+1/2-1/p}|v|_{W^{k+1,p}(\Omega)}. 
	\end{equation}
	For the normal derivative, the inverse trace estimate~\eqref{eq:grad-inverse-trace} applied face by face combined with~\eqref{eq:zh-bdry-layer} gives $$\|\partial_nz_h\|_{L^2(\Gamma)} \le Ch^{-1/2}\|\nabla z_h\|_{L^2(\Omega_h^\Gamma)} \le Ch^{k-1/p}|v|_{W^{k+1,p}(\Omega)}.$$ Using the triangle inequality for $\partial_n\rho = \partial_n\eta_0 - \partial_n\psi_E - \partial_n z_h$, along with the bounds for $\partial_n\eta_0$ from~\eqref{eq:interp-tang} and $\partial_n\psi_E$ from~\eqref{eq:psiE-bounds}, we obtain
	\[
	h^{1/2}\|\partial_n\rho\|_{L^2(\Gamma)} \le Ch^{k+1/2-1/p}|v|_{W^{k+1,p}(\Omega)}. 
	\]
	Combining this with~\eqref{eq:rho-bdry-layer-grad} and~\eqref{eq:rho-trace} completes the proof of~\eqref{eq:Ritz-boundary-layer}.

\end{proof}

\begin{remark}
	A direct route via a $W^{1,p}$ stability estimate for the nonhomogeneous Ritz projection $\mathcal R_h$ applied to $\eta=v-J_hv$ would, at $p=\infty$, carry a $|\log h|$ factor. The harmonic lift reduction $z=\eta-q\in H_0^1(\Omega)$ avoids this by passing to the homogeneous projection $\mathcal R_h^0$, on which the log-free pointwise gradient estimate of Diening--Rolfes--Salgado~\cite[Theorem~3.1]{DieningRolfesSalgado2024} applies, and the $L^p$ boundedness of the Hardy--Littlewood maximal operator absorbs the boundary-layer factor $h^{1/2-1/p}$ uniformly for $p\in[2,\infty]$.
\end{remark}

\subsection{Boundary-strip error estimate}

We transfer the boundary-layer Ritz estimate of \Cref{lem:nh-ritz-boundary-layer} to the discrete error $u-u_h$ via the energy norm.

\begin{lemma}\label{lem:nitsche-energy-boundary}
	Under Assumption~\ref{assump:standing}, let $\alpha\ge 1$. Let $u\in W^{k+1,p}(\Omega)$ for some $p\in[2,\infty]$ be the exact solution of~\eqref{org}, and let $u_h\in V_h$ solve~\eqref{eq:prelim-nitsche}.  Then the bound~\eqref{eq:critical-strip-error-opening} holds:
	\begin{equation*}
		\begin{split}
			&\|\nabla(u-u_h)\|_{L^2(\Omega_h^\Gamma)} + h^{1/2}\|\partial_n(u-u_h)\|_{L^2(\Gamma)} + h^{-1/2}\|u-u_h\|_{L^2(\Gamma)} \\
			&\qquad \le Ch^{k+1/2-1/p}|u|_{W^{k+1,p}(\Omega)}. 
		\end{split}
	\end{equation*}
	Moreover, the projected-boundary Ritz error satisfies the sharper estimate
	\begin{equation}\label{eq:critical-Ritz-boundary}
		\|\mathcal R_hu-u_h\|_{L^2(\Gamma)}\le Ch^{k+(1+\alpha)/2-1/p}|u|_{W^{k+1,p}(\Omega)}. 
	\end{equation}
\end{lemma}
\begin{proof}
	{
		Decompose $u-u_h = (u-\mathcal R_hu) + (\mathcal R_hu-u_h) := \rho + r_h$. We claim that both $\rho$ and $r_h$ satisfy the target bound:
		\begin{equation}\label{eq:summand-target}
			\|\nabla w\|_{L^2(\Omega_h^\Gamma)} + h^{1/2}\|\partial_n w\|_{L^2(\Gamma)} + h^{-1/2}\|w\|_{L^2(\Gamma)} \le Ch^{k+1/2-1/p}|u|_{W^{k+1,p}(\Omega)}.
		\end{equation}
		The triangle inequality then yields~\eqref{eq:critical-strip-error-opening}
	}

	For {$\rho=u-\mathcal R_hu$}, the bound~\eqref{eq:summand-target} is precisely~\eqref{eq:Ritz-boundary-layer} applied with $v=u$ (the same three quantities, written in different order).
	
	For {$r_h = \mathcal R_hu-u_h\in V_h$}, we use the energy norm together with the trace bounds
	\begin{equation}\label{eq:vh-trace-bounds}
		\|v_h\|_{L^2(\Gamma)}\le h^{\alpha/2}\vertiii{v_h}_h,\quad
		\|\partial_n v_h\|_{L^2(\Gamma)}\le Ch^{-1/2}\vertiii{v_h}_h \quad\text{for } v_h\in V_h,
	\end{equation}
	valid by~\eqref{eq:reference-nitsche-norm} (with $\alpha_\ast=\alpha$ for $\alpha\ge 1$) and the inverse trace estimate~\eqref{eq:grad-inverse-trace} applied with $\|\nabla v_h\|_{L^2(\Omega)}\le\vertiii{v_h}_h$. Galerkin orthogonality~\eqref{eq:galerkin-orthogonality} with {$v_h =r_h=\mathcal R_hu-u_h$} and the diagonal identity~\eqref{eq:diag-cancel} give
	\[
	{\vertiii{r_h}_h^2=\vertiii{\mathcal R_hu-u_h}_h^2 = -a(u-\mathcal R_hu,\,\mathcal R_hu-u_h)=-a(\rho,r_h).} 
	\]

	{
		We now bound the bilinear form $a(\rho, r_h)$ by expanding it via~\eqref{eq:ah-def}:
		\begin{equation*}
			a(\rho,r_h):=(\nabla\rho,\nabla r_h)_\Omega-\langle\partial_n \rho,r_h\rangle_\Gamma
			+\langle \rho,\partial_n r_h\rangle_\Gamma+h^{-\alpha}\langle \rho,r_h\rangle_\Gamma,
		\end{equation*}
		\begin{itemize}
			\item \textbf{$\langle \rho,r_h\rangle_\Gamma$:} Vanishes by boundary orthogonality~\eqref{eq:Ritz-boundary-orthogonality}.
			\item \textbf{$(\nabla\rho,\nabla r_h)$:} Let $\mu_h = r_h|_\Gamma$ and define $w_h := r_h - E_h(\mu_h) \in V_h^0$. The Ritz orthogonality $(\nabla\rho, \nabla w_h)_\Omega = 0$ allows us to rewrite the volume integral as $(\nabla\rho, \nabla E_h(\mu_h))_\Omega$. Since $E_h(\mu_h)$ is supported in $\overline{\Omega_h^\Gamma}$ by~\eqref{eq:disc-lift-trace}, we can apply Cauchy--Schwarz, the boundary layer estimate~\eqref{eq:Ritz-boundary-layer} for $\rho$, and the discrete lift bound~\eqref{eq:disc-lift-W12} to bound this term.
			\item \textbf{$\langle\partial_n \rho,r_h\rangle_\Gamma
				$ and $\langle \rho,\partial_n r_h\rangle_\Gamma$:} Bounded directly via Cauchy--Schwarz, applying~\eqref{eq:Ritz-boundary-layer} for $\rho$ and the trace bounds~\eqref{eq:vh-trace-bounds} for $r_h$.
		\end{itemize}
	}
	Using $h^{\alpha/2-1/2}\le 1$ for $\alpha\ge 1$ to consolidate the volume and first-boundary terms, both initially scaling as $h^{k+\alpha/2-1/p}$ and 
	{combining these estimates, we obtain
		\[
		\vertiii{r_h}_h^2=|a(\rho,\, r_h)| \le Ch^{k+1/2-1/p}|u|_{W^{k+1,p}(\Omega)}\vertiii{r_h}_h.
		\]
		Dividing by $\vertiii{r_h}_h$ (trivial if zero) yields the energy bound:
		\begin{equation}\label{eq:eh-energy-bound}
			\vertiii{r_h}_h \le Ch^{k+1/2-1/p}|u|_{W^{k+1,p}(\Omega)}. 
		\end{equation}
	}
	
	Combining this with~\eqref{eq:vh-trace-bounds} and $\|\nabla v_h\|_{L^2(\Omega)}\le\vertiii{v_h}_h$, the three terms of~\eqref{eq:summand-target} for $r_h = \mathcal R_hu-u_h$ are each bounded by $Ch^{k+1/2-1/p}|u|_{W^{k+1,p}(\Omega)}$ (using $h^{(\alpha-1)/2}\le 1$ for $\alpha\ge 1$ in the $h^{-1/2}\|\cdot\|_{L^2(\Gamma)}$ term). The sharper bound~\eqref{eq:critical-Ritz-boundary} follows by not absorbing this $h^{(\alpha-1)/2}$ factor:
	\[
	{
		\|r_h\|_{L^2(\Gamma)} \le h^{\alpha/2}\vertiii{r_h}_h \le Ch^{k+(1+\alpha)/2-1/p}|u|_{W^{k+1,p}(\Omega)}. }
	\]
\end{proof}

\subsection{Sharp estimate of \texorpdfstring{$T_3$}{T3} for \texorpdfstring{$\alpha\ge 1$}{α≥1}}

Applying the boundary-functional identity~\eqref{eq:boundary-functional-identity} with the boundary-supported lifting (\Cref{lem:disc-lift}) and the dual-lift scales (\Cref{cor:dual-lift-scales}), we decompose $T_3$ into three pieces and bound each using \Cref{lem:nitsche-energy-boundary}.

\begin{lemma}\label{lem:T3-sharp}
	Under Assumption~\ref{assump:standing}, let $\alpha\ge 1$. Let $u\in W^{k+1,p}(\Omega)$ for some $p\in[2,\infty]$ be the exact solution of~\eqref{org}, and let $u_h\in V_h$ solve~\eqref{eq:prelim-nitsche}.  Then the projected boundary term $T_3$ from \Cref{lem:standard-duality-T3}, taken in its Ritz form~\eqref{eq:T3-Ritz-form}, satisfies
	\begin{equation}\label{eq:T3-sharp}
		|T_3| \le C h^{k+\alpha-1/p}|u|_{W^{k+1,p}(\Omega)}\|u-u_h\|_{L^2(\Omega)}. 
	\end{equation}
\end{lemma}

\begin{proof}
	Let $\Psi\in H^2(\Omega)\cap H_0^1(\Omega)$ solve $-\Delta\Psi=u-u_h$ in $\Omega$, and set
	\[
	v_h := E_h(\Pi_h^\Gamma\partial_n\Psi)\in V_h. 
	\]
	By~\eqref{eq:disc-lift-trace}, $v_h|_\Gamma = \Pi_h^\Gamma\partial_n\Psi\in M_h$, so $v_h$ is admissible in the boundary-functional identity~\eqref{eq:boundary-functional-identity}, and $\operatorname{supp}(v_h)\subset\overline{\Omega_h^\Gamma}$. Using $\langle\partial_n(u-u_h),v_h\rangle_\Gamma = \langle\partial_n(u-u_h),\Pi_h^\Gamma\partial_n\Psi\rangle_\Gamma$, the identity decomposes $T_3 = T_{31}+T_{32}+T_{33}$ with
	\begin{align*}
		T_{31} &:= -h^\alpha(\nabla(u-u_h),\nabla v_h)_\Omega, \\
		T_{32} &:= h^\alpha\langle\partial_n(u-u_h),\Pi_h^\Gamma\partial_n\Psi\rangle_\Gamma, \\
		T_{33} &:= -h^\alpha\langle u-u_h,\partial_n v_h\rangle_\Gamma. 
	\end{align*}
	Applying Cauchy--Schwarz, with the a.e.\ vanishing of $v_h$ on $\Omega\setminus\Omega_h^\Gamma$ from~\eqref{eq:disc-lift-trace} restricting the $T_{31}$ integral to $\Omega_h^\Gamma$,
	\begin{subequations}\label{eq:T3j-CS}
		\begin{align}
			|T_{31}| &\le h^\alpha\,\|\nabla(u-u_h)\|_{L^2(\Omega_h^\Gamma)}\,\|\nabla v_h\|_{L^2(\Omega)}, \\
			|T_{32}| &\le h^\alpha\,\|\partial_n(u-u_h)\|_{L^2(\Gamma)}\,\|\Pi_h^\Gamma\partial_n\Psi\|_{L^2(\Gamma)}, \\
			|T_{33}| &\le h^\alpha\,\|u-u_h\|_{L^2(\Gamma)}\,\|\partial_n v_h\|_{L^2(\Gamma)}. 
		\end{align}
	\end{subequations}
	The bound~\eqref{eq:critical-strip-error-opening} gives the three Nitsche-error factors:
	\begin{align*}
		\|\nabla(u-u_h)\|_{L^2(\Omega_h^\Gamma)} &\le Ch^{k+1/2-1/p}|u|_{W^{k+1,p}(\Omega)}, \\
		\|\partial_n(u-u_h)\|_{L^2(\Gamma)} &\le Ch^{k-1/p}|u|_{W^{k+1,p}(\Omega)}, \\
		\|u-u_h\|_{L^2(\Gamma)} &\le Ch^{k+1-1/p}|u|_{W^{k+1,p}(\Omega)},
	\end{align*}
	the second after dividing by $h^{1/2}$ and the third after multiplying by $h^{1/2}$. \Cref{cor:dual-lift-scales}, via~\eqref{eq:dual-lift-W12},~\eqref{eq:Pi-trace-bound}, and~\eqref{eq:dual-lift-normal}, gives the three dual-lift factors:
	\begin{align*}
		\|\nabla v_h\|_{L^2(\Omega)} &\le Ch^{-1/2}\|\Psi\|_{H^2(\Omega)}, \\
		\|\Pi_h^\Gamma\partial_n\Psi\|_{L^2(\Gamma)} &\le C\|\Psi\|_{H^2(\Omega)}, \\
		\|\partial_n v_h\|_{L^2(\Gamma)} &\le Ch^{-1}\|\Psi\|_{H^2(\Omega)}. 
	\end{align*}
	Each corresponding product equals $Ch^{k+\alpha-1/p}|u|_{W^{k+1,p}(\Omega)}\|\Psi\|_{H^2(\Omega)}$, so substituting into~\eqref{eq:T3j-CS} and summing,
	\[
	|T_3|\le|T_{31}|+|T_{32}|+|T_{33}|\le Ch^{k+\alpha-1/p}|u|_{W^{k+1,p}(\Omega)}\|\Psi\|_{H^2(\Omega)}. 
	\]
	The elliptic regularity~\eqref{eq:standard-dual-regularity}, $\|\Psi\|_{H^2(\Omega)}\le C\|u-u_h\|_{L^2(\Omega)}$, proves~\eqref{eq:T3-sharp}.
\end{proof}

\begin{remark}
	The Galerkin orthogonality argument above is essential only for $\alpha>1$. At the borderline $\alpha=1$, the sharper Ritz boundary estimate~\eqref{eq:critical-Ritz-boundary} reduces to $\|\mathcal R_hu-u_h\|_{L^2(\Gamma)}\le Ch^{k+1-1/p}|u|_{W^{k+1,p}(\Omega)}$, and combining this with Cauchy--Schwarz on $T_3$ in its Ritz form~\eqref{eq:T3-Ritz-form}, the $L^2(\Gamma)$ stability of $\Pi_h^\Gamma$, the trace theorem~\eqref{trace_theorem}, and the elliptic regularity~\eqref{eq:standard-dual-regularity} already gives
	\[
	|T_3|\le\|\mathcal R_hu-u_h\|_{L^2(\Gamma)}\|\Pi_h^\Gamma\partial_n\Psi\|_{L^2(\Gamma)}\le Ch^{k+1-1/p}|u|_{W^{k+1,p}(\Omega)}\|u-u_h\|_{L^2(\Omega)},
	\]
	which agrees with~\eqref{eq:T3-sharp} at $\alpha=1$. For $\alpha>1$, this direct route yields only $h^{k+(1+\alpha)/2-1/p}$, missing $(\alpha-1)/2$ orders: recovering the full penalty factor $h^\alpha$ requires extracting $T_3$ from Galerkin orthogonality via the boundary-functional identity~\eqref{eq:boundary-functional-identity}, as done in the proof.
\end{remark}

\begin{remark}
	The difference between Sections~\ref{sec:traditional-duality} and~\ref{sec:main} is the treatment of $T_3$. The traditional argument applies Cauchy--Schwarz directly to the boundary inner product defining $T_3$, extracting only the square root $h^{\alpha/2}$ from the boundary energy norm. The corrected argument uses the boundary-functional identity~\eqref{eq:boundary-functional-identity} to express $T_3$ as $h^\alpha$ times a sum of volume and trace residuals, exploiting that the penalty multiplier $h^{-\alpha}$ acting on the boundary inner product yields the full factor $h^\alpha$. The strip estimate~\eqref{eq:critical-strip-error-opening} then supplies the additional factor $h^{1/2-1/p}$ via $W^{k+1,p}$ regularity.
\end{remark}

\subsection{Sharp \texorpdfstring{$L^2$}{L2} estimate for \texorpdfstring{$\alpha\ge 1$}{α≥1}}

The sharp $T_3$ bound of \Cref{lem:T3-sharp}, combined with the bounds on $T_1$ and $T_2$ already established in \Cref{lem:standard-duality-T3}, immediately yields the main result of this section.

\begin{theorem}\label{thm:main-L2}
	Under Assumption~\ref{assump:standing}, let $\alpha\ge 1$. Let $u\in W^{k+1,p}(\Omega)$ for some $p\in[2,\infty]$ be the exact solution of~\eqref{org}, and let $u_h\in V_h$ solve~\eqref{eq:prelim-nitsche}.  Then
	\begin{equation}\label{eq:main-bound-alpha}
		\|u-u_h\|_{L^2(\Omega)}
		\le C h^{\min\{k+1,\,k+\alpha-1/p\}} |u|_{W^{k+1,p}(\Omega)}. 
	\end{equation}
\end{theorem}

\begin{proof}
	By the dual identity~\eqref{eq:standard-duality-T123}, the universal bound~\eqref{eq:T1-T2-universal} combined with the H\"older embedding $|u|_{H^{k+1}(\Omega)}\le C|u|_{W^{k+1,p}(\Omega)}$ for $p\ge 2$ and the elliptic regularity~\eqref{eq:standard-dual-regularity}, and the sharp $T_3$ estimate~\eqref{eq:T3-sharp} from \Cref{lem:T3-sharp},
	\[
	\|u-u_h\|_{L^2(\Omega)}^2
	\le C\bigl(h^{k+1}+h^{k+\alpha-1/p}\bigr)|u|_{W^{k+1,p}(\Omega)}\|u-u_h\|_{L^2(\Omega)}. 
	\]
	Dividing by $\|u-u_h\|_{L^2(\Omega)}$ proves~\eqref{eq:main-bound-alpha}.
\end{proof}

\begin{remark}
	For the classical stabilized scaling $\alpha=1$ used in the original statement and in the numerical experiments, \Cref{thm:main-L2} specializes to
	\begin{equation}
		\|u-u_h\|_{L^2(\Omega)}
		\le C h^{k+1-1/p}|u|_{W^{k+1,p}(\Omega)}. 
	\end{equation}
\end{remark}

\section{Sharp \texorpdfstring{$L^2$}{L2} estimate for all \texorpdfstring{$\alpha$}{α}}\label{sec:weak-penalty}

\Cref{thm:main-L2} has already established the sharp $L^2$ estimate for $\alpha\ge 1$.  It remains to treat the weak-penalty range
\[
\alpha\in\{-\infty\}\cup(-\infty,1),
\]
where $\alpha=-\infty$ denotes the penalty-free method; combined with~\Cref{thm:main-L2}, the result of this section will complete the unified statement~\Cref{thm:unified-L2}.  The dual identity of \Cref{lem:standard-duality-T3} again isolates the projected boundary term in its Ritz form~\eqref{eq:T3-Ritz-form},
\[
T_3=\langle\mathcal R_hu-u_h,\Pi_h^\Gamma\partial_n\Psi\rangle_\Gamma. 
\]

The strategy, however, is simpler than in Section~\ref{sec:main}.  There, the penalty factor $h^\alpha$ had to be extracted by routing $T_3$ through Galerkin orthogonality; here no such gain exists, since the penalty is absent at $\alpha=-\infty$ and weaker than the reference boundary scaling for finite $\alpha<1$.  A direct Cauchy--Schwarz on $T_3$ then reduces the entire analysis to sharpening the projected-boundary Ritz error $\|\mathcal R_hu-u_h\|_{L^2(\Gamma)}$.

The direct bound from \Cref{lem:traditional-energy} is only $Ch^{k+1/2}|u|_{H^{k+1}(\Omega)}$, and upgrading the regularity globally to $W^{k+1,p}$ contributes nothing beyond a mesh-independent constant: as in Section~\ref{sec:main}, the higher regularity must be harnessed locally rather than globally.  The same boundary-strip mechanism applies here, now turned on the Ritz energy norm $\vertiii{\mathcal R_hu-u_h}_h$ via the inf--sup condition: H\"older's inequality localized to the strip $\Omega_h^\Gamma$ supplies a factor $h^{1/2-1/p}$.

Section~\ref{sec:weak-bdry-est} carries out the argument; Section~\ref{sec:weak-L2-est} combines it with the dual identity to produce the $L^2$ rate.

\subsection{Sharp estimate of \texorpdfstring{$T_3$}{T3} for \texorpdfstring{$\alpha<1$}{α<1}}\label{sec:weak-bdry-est}

\begin{lemma}\label{lem:weak-boundary-estimate}
	Under Assumption~\ref{assump:standing}, let $\alpha\in\{-\infty\}\cup(-\infty,1)$ and use meshes with $h<h_0$ as in \Cref{thm:weak-infsup}. Let $u\in W^{k+1,p}(\Omega)$ for some $p\in[2,\infty]$ be the exact solution of~\eqref{org}, and let $u_h\in V_h$ solve~\eqref{eq:prelim-nitsche}. Then
	\begin{equation}\label{eq:weak-boundary-est}
		\|\mathcal R_hu-u_h\|_{L^2(\Gamma)}\le C h^{k+1-1/p}|u|_{W^{k+1,p}(\Omega)}. 
	\end{equation}
	Consequently, the projected boundary term $T_3$ from \Cref{lem:standard-duality-T3} satisfies
	\begin{equation}\label{eq:weak-T3-bound}
		|T_3|\le C h^{k+1-1/p}|u|_{W^{k+1,p}(\Omega)}\|u-u_h\|_{L^2(\Omega)}. 
	\end{equation}
\end{lemma}

\begin{proof}
	We compare the discrete solution with the projected-boundary Ritz projection $\mathcal R_hu$ defined in~\eqref{eq:nh-ritz-def}.  Since $\mathcal R_hu-u_h\in V_h$, the inf--sup condition and Galerkin orthogonality give
	\[
	\vertiii{\mathcal R_hu-u_h}_h
	\le C\sup_{0\ne w_h\in V_h}
	\frac{|a(u-\mathcal R_hu,w_h)|}{\vertiii{w_h}_h}. 
	\]
	For $w_h\in V_h$, the boundary orthogonality $\langle u-\mathcal R_hu,w_h\rangle_\Gamma=0$ (valid since $w_h|_\Gamma\in M_h$) makes the penalty part of $a(u-\mathcal R_hu,w_h)$ vanish for every weak penalty, including the penalty-free endpoint where the term is absent.  Let $E_h(w_h|_\Gamma)$ be the boundary lifting of \Cref{lem:disc-lift}.  Since $w_h-E_h(w_h|_\Gamma)\in V_h^0$, the Ritz orthogonality of $u-\mathcal R_hu$ recasts the volume contribution as
	\[
	(\nabla(u-\mathcal R_hu),\nabla w_h)_\Omega
	=(\nabla(u-\mathcal R_hu),\nabla E_h(w_h|_\Gamma))_\Omega,
	\]
	supported in $\overline{\Omega_h^\Gamma}$.  Set $\rho:=u-\mathcal R_hu$ and $\mu_h:=w_h|_\Gamma$.  The preceding identities give the exact decomposition
	\begin{equation*}
		a(\rho,w_h)
		=(\nabla\rho,\nabla E_h\mu_h)_\Omega
		-\langle\partial_n\rho,w_h\rangle_\Gamma
		+\langle\rho,\partial_n w_h\rangle_\Gamma .
	\end{equation*}
	We estimate the three terms separately.  Put
	$A_u:=h^{k+1/2-1/p}|u|_{W^{k+1,p}(\Omega)}$.  By~\eqref{eq:Ritz-boundary-layer} and~\eqref{eq:disc-lift-W12},
	\begin{align*}
		|(\nabla\rho,\nabla E_h\mu_h)_\Omega|
		&\le \|\nabla\rho\|_{L^2(\Omega_h^\Gamma)}\,\|\nabla E_h\mu_h\|_{L^2(\Omega)} \\
		&\le C A_u\, h^{-1/2}\|w_h\|_{L^2(\Gamma)}
		\le C A_u\vertiii{w_h}_h .
	\end{align*}
	Again by~\eqref{eq:Ritz-boundary-layer} and by the reference trace control in $\vertiii{\cdot}_h$ for $\alpha<1$,
	\begin{align*}
		|\langle\partial_n\rho,w_h\rangle_\Gamma|
		&\le \|\partial_n\rho\|_{L^2(\Gamma)}\|w_h\|_{L^2(\Gamma)} \\
		&\le C h^{-1/2}A_u\,h^{1/2}\vertiii{w_h}_h
		\le C A_u\vertiii{w_h}_h .
	\end{align*}
	Finally, the trace part of~\eqref{eq:Ritz-boundary-layer} and the inverse trace estimate~\eqref{eq:grad-inverse-trace} give
	\begin{align*}
		|\langle\rho,\partial_n w_h\rangle_\Gamma|
		&\le \|\rho\|_{L^2(\Gamma)}\|\partial_n w_h\|_{L^2(\Gamma)} \\
		&\le C h^{1/2}A_u\,h^{-1/2}\|\nabla w_h\|_{L^2(\Omega)}
		\le C A_u\vertiii{w_h}_h .
	\end{align*}
	Thus $|a(u-\mathcal R_hu,w_h)|\le C h^{k+1/2-1/p}|u|_{W^{k+1,p}(\Omega)}\vertiii{w_h}_h$, and substitution into the inf--sup bound yields the energy estimate
	\[
	\vertiii{\mathcal R_hu-u_h}_h\le Ch^{k+1/2-1/p}|u|_{W^{k+1,p}(\Omega)}. 
	\]
	The reference trace control then gives~\eqref{eq:weak-boundary-est}:
	\[
	\|\mathcal R_hu-u_h\|_{L^2(\Gamma)}\le h^{1/2}\vertiii{\mathcal R_hu-u_h}_h\le Ch^{k+1-1/p}|u|_{W^{k+1,p}(\Omega)},
	\]
	using $\alpha_\ast=1$ for $\alpha<1$.
	
	For the consequence~\eqref{eq:weak-T3-bound}, let $\Psi\in H^2(\Omega)\cap H_0^1(\Omega)$ be the dual solution from \Cref{lem:standard-duality-T3}.  Using the Ritz form~\eqref{eq:T3-Ritz-form}, Cauchy--Schwarz, the $L^2(\Gamma)$ stability of $\Pi_h^\Gamma$, the trace theorem~\eqref{trace_theorem}, the elliptic regularity~\eqref{eq:standard-dual-regularity}, and~\eqref{eq:weak-boundary-est} give
	\[
	|T_3|\le \|\mathcal R_hu-u_h\|_{L^2(\Gamma)}\|\Pi_h^\Gamma\partial_n\Psi\|_{L^2(\Gamma)}
	\le Ch^{k+1-1/p}|u|_{W^{k+1,p}(\Omega)}\|u-u_h\|_{L^2(\Omega)}. 
	\]
\end{proof}

\subsection{Sharp \texorpdfstring{$L^2$}{L2} estimate for \texorpdfstring{$\alpha<1$}{α<1}}\label{sec:weak-L2-est}

\begin{theorem}\label{thm:weak-L2}
	Under Assumption~\ref{assump:standing}, let $\alpha\in\{-\infty\}\cup(-\infty,1)$ and use meshes with $h<h_0$ as in \Cref{thm:weak-infsup}. Let $u\in W^{k+1,p}(\Omega)$ for some $p\in[2,\infty]$ be the exact solution of~\eqref{org}, and let $u_h\in V_h$ solve~\eqref{eq:prelim-nitsche}. Then
	\begin{equation}\label{eq:weak-L2-rate}
		\|u-u_h\|_{L^2(\Omega)}
		\le C h^{k+1-1/p}|u|_{W^{k+1,p}(\Omega)}. 
	\end{equation}
\end{theorem}

\begin{proof}
	By the dual identity~\eqref{eq:standard-duality-T123}, the universal bound~\eqref{eq:T1-T2-universal} combined with the H\"older embedding $|u|_{H^{k+1}(\Omega)}\le C|u|_{W^{k+1,p}(\Omega)}$ for $p\ge 2$ and the elliptic regularity~\eqref{eq:standard-dual-regularity}, and the $T_3$ estimate~\eqref{eq:weak-T3-bound} from \Cref{lem:weak-boundary-estimate},
	\[
	\|u-u_h\|_{L^2(\Omega)}^2
	\le C h^{k+1-1/p}|u|_{W^{k+1,p}(\Omega)}\|u-u_h\|_{L^2(\Omega)}. 
	\]
	Dividing by $\|u-u_h\|_{L^2(\Omega)}$ proves~\eqref{eq:weak-L2-rate}.
\end{proof}

\subsection{Sharp \texorpdfstring{$L^2$}{L2} estimate for all \texorpdfstring{$\alpha$}{α}}

Combining \Cref{thm:main-L2} for $\alpha\ge1$ with \Cref{thm:weak-L2} for $\alpha<1$ yields the unified $L^2$ estimate across the entire penalty scale, the main result of the paper.

\begin{theorem}[Unified $L^2$ estimate]\label{thm:unified-L2}
	Under Assumption~\ref{assump:standing}, let $\alpha\in\{-\infty\}\cup\mathbb R$; when $\alpha<1$, use meshes with $h<h_0$ as in \Cref{thm:weak-infsup}. Let $u\in W^{k+1,p}(\Omega)$ for some $p\in[2,\infty]$ be the exact solution of~\eqref{org}, and let $u_h\in V_h$ solve~\eqref{eq:prelim-nitsche}. Then
	\begin{equation}\label{eq:full-penalty-scale-rate}
		\begin{split}
			\|u-u_h\|_{L^2(\Omega)} \le C h^{r}|u|_{W^{k+1,p}(\Omega)}, \quad 
			r = \min\bigl\{k+1,\,k+\max\{1,\alpha\}-1/p\bigr\}.
		\end{split}
	\end{equation}
\end{theorem}

\begin{figure}[ht]
	\centering
	\begin{minipage}[t]{0.55\linewidth}
		\vspace{0pt}
		\centering
		\includegraphics[width=\linewidth]{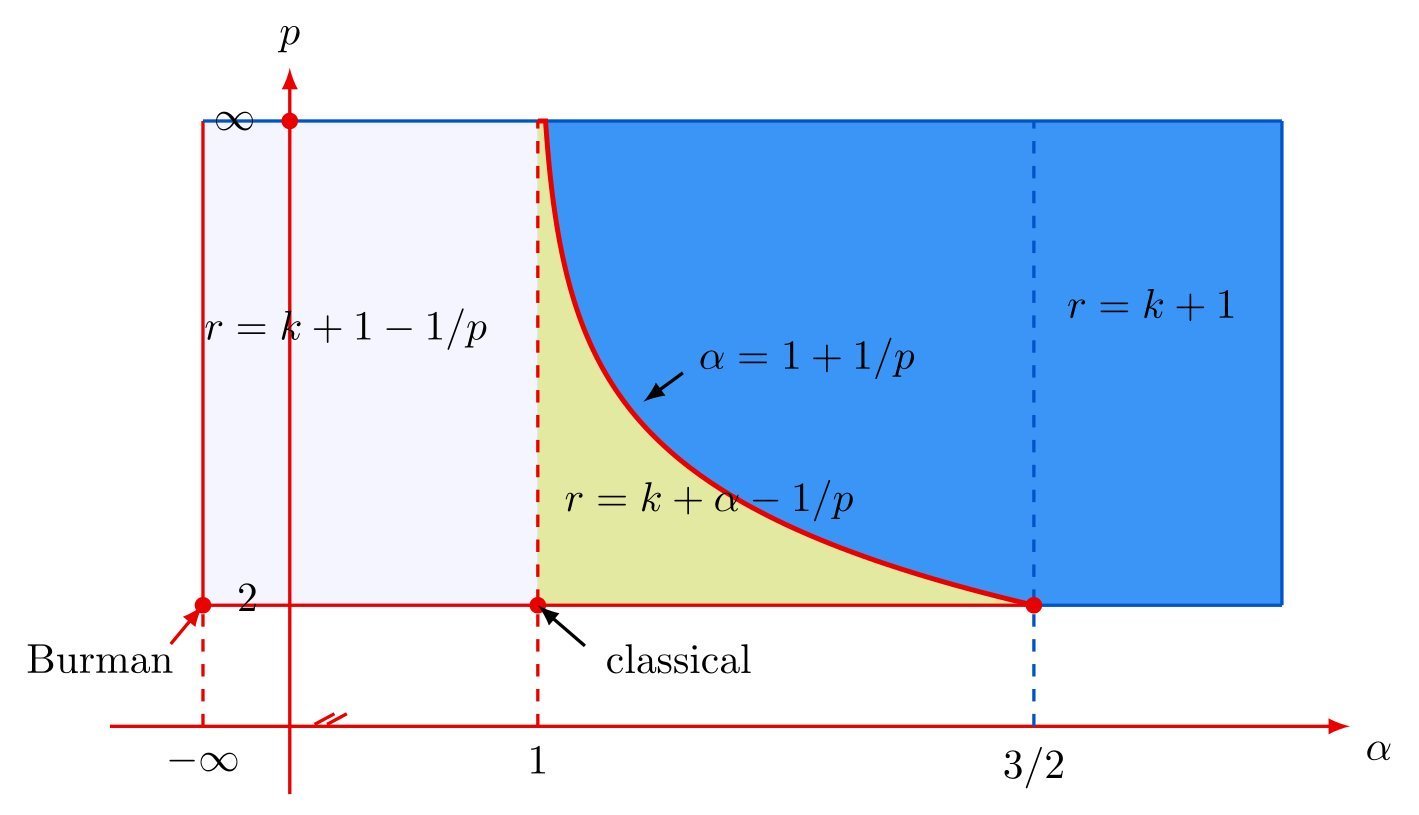}
		\caption{Coverage diagram for \Cref{thm:unified-L2} in the $(\alpha,p)$-plane.}
		\label{fig:coverage-diagram}
	\end{minipage}\hfill
	\begin{minipage}[t]{0.42\linewidth}
		\vspace{0pt}
		\centering
		\captionof{table}{Penalty regimes and $L^2$ convergence orders from \Cref{thm:unified-L2}.}
		\label{tab:penalty-summary}
		\vspace{4pt}
		\renewcommand{\arraystretch}{1.6}
		\begin{tabular}{c|c}
			\hline
			\textbf{Penalty regime} & \textbf{Order $r$} \\
			\hline\hline
			$\alpha\le 1$ & $k+1-1/p$ \\
			$1<\alpha<1+1/p$ & $k+\alpha-1/p$ \\
			$\alpha\ge 1+1/p$ & $k+1$ \\
			\hline
		\end{tabular}
	\end{minipage}
\end{figure}

In particular, the optimal rate $r=k+1$ is attained if and only if $\alpha\ge 1+1/p$ or $p=\infty$. The rate  $r$ across the penalty scale is summarized in \Cref{fig:coverage-diagram} and \Cref{tab:penalty-summary}.

\section{Numerical experiments}\label{sec:numerics}

For piecewise linear elements, \Cref{thm:unified-L2} predicts the unified $L^2$ rate
\begin{align}\label{eq:numerics-predicted-rate}
	\|u-u_h\|_{L^2(\Omega)} \le C\, h^{r}\,|u|_{W^{2,p}(\Omega)},
	\qquad r = \min\bigl\{2,\,1+\max\{1,\alpha\}-1/p\bigr\}.
\end{align}

\paragraph{Setup} All computations are carried out in C++ on top of the MFEM finite element library~\cite{anderson2021mfem}. We discretize the Poisson problem~\eqref{org} on $\Omega=(0,1)^d$, in both $d=2$ and $d=3$, with continuous piecewise linear Lagrange elements ($k=1$), so that the prediction~\eqref{eq:numerics-predicted-rate} is the second-order target tested throughout. The resulting linear systems are solved by GMRES~\cite{MR848568} for the nonsymmetric formulations and by the conjugate gradient method~\cite{hestenes1952methods} for the symmetric reference, with BoomerAMG~\cite{MR972756} from the \texttt{hypre} library~\cite{MR2267935} as preconditioner in both cases. The penalty coefficient is $\gamma=h^{-\alpha}$ throughout, reducing to the classical Nitsche scaling $\gamma=1/h$ at $\alpha=1$ and vanishing identically at Burman's penalty-free endpoint $\alpha=-\infty$; the adjoint-consistent symmetric variant $\beta=-1$, included as a reference, is run with $\gamma=10^6/h$. Successive uniform refinements bring the discrete problem up to roughly $2.7\times 10^8$ degrees of freedom in two dimensions and $1.7\times 10^7$ in three dimensions. In \Cref{tab:numerics-2D,tab:numerics-3D}, the column labeled ``pred'' records the predicted rate from~\eqref{eq:numerics-predicted-rate}, and each subsequent column reports the relative $L^2$ error on the top row, together with the observed convergence rate against the immediately preceding refinement level on the bottom row.

\paragraph{The test family} Define
\begin{equation}\label{eq:numerics-family}
	u_p(\mathbf x) := \big[x_1(1-x_1)\big]^{\,2-1/p}\cdot\prod_{i=2}^{d} x_i^{\,2-1/p},\qquad p\in[2,\infty],\quad d\in\{2,3\}.
\end{equation}
The omission of \((1-x_i)\) for \(i=2,\ldots,d\) leaves \(u_p\)
nonzero on the faces \(\{x_i=1\}\), producing genuinely nonhomogeneous
Dirichlet data; at \(p=\infty\), \(u_\infty\) is a polynomial in
\(W^{2,\infty}\). A direct calculation shows that, for
\(2\le p<\infty\), the singular factor \(x_i^{2-1/p}\) has second
derivative behaving like \(x_i^{-1/p}\) near \(x_i=0\), so, for every
sufficiently small \(\varepsilon>0\),
\[
u_p\in W^{2,p-\varepsilon}(\Omega),\qquad
u_p\notin W^{2,p}(\Omega).
\]
Thus \(p\) is the limiting integrability index of \(u_p\) for finite
\(p\), while at the endpoint \(p=\infty\), \(u_\infty\) is a polynomial
and lies in \(W^{2,\infty}(\Omega)\). For finite \(p>2\), we apply
\Cref{thm:unified-L2} with \(p_\varepsilon:=p-\varepsilon\), where
\(\varepsilon>0\) is chosen so small that \(p_\varepsilon\ge2\), and
then report the limiting predicted rate obtained as
\(\varepsilon\downarrow0\). For \(p=\infty\), \Cref{thm:unified-L2} applies
directly. Hence the predicted rate \eqref{eq:numerics-predicted-rate} reduces to \(2-1/p\) at
\(\alpha\in\{-\infty,1\}\) and to
\(\min\{2,1+\alpha-1/p\}\) at \(1<\alpha\le1.5\).

\subsection{Sharpness of \texorpdfstring{$r=2-1/p$}{r=2-1/p} at \texorpdfstring{$\alpha\in\{1,-\infty\}$}{α∈\{1,-∞\}}}\label{sec:numerics-endpoints}

\Cref{tab:numerics-2D,tab:numerics-3D} report convergence histories at $p\in\{2.0004,\,4,\,\infty\}$ for both nonsymmetric endpoints $\alpha=1$ and $\alpha=-\infty$ (which share $\max\{1,\alpha\}=1$), with the symmetric method as adjoint-consistent reference, in $d\in\{2,3\}$. The data confirm \Cref{thm:unified-L2}'s prediction $r=2-1/p$ entrywise: the half-order loss $h^{3/2}$ at the rough endpoint $p=2.0004$, where $u_p$ has essentially only $H^2$ regularity; the optimal rate $h^2$ at the smooth endpoint $p=\infty$, where $u_\infty\in W^{2,\infty}$; and continuous interpolation between them at the intermediate sample $p=4$. Together with the symmetric reference, which attains $r=2$ throughout the same data, this resolves the open question: under merely $H^{k+1}$ regularity the half-order loss is essential and intrinsic to the nonsymmetric method, while the optimal rate consistently observed on smooth test problems is explained by their full $W^{k+1,\infty}$ regularity, not by a limitation of the standard analysis.

\begin{table}[htbp]
	\centering
	\caption{2D convergence histories on the family~\eqref{eq:numerics-family}: relative $L^2$ error (top) and observed rate (bottom). The ``pred'' column is the rate from~\eqref{eq:numerics-predicted-rate}.}
	\label{tab:numerics-2D}
	\resizebox{\textwidth}{!}{%
		\begin{tabular}{c|c|c|ccccc}
			\hline
			\multirow{2}{*}{$p$} & \multirow{2}{*}{Method} & \multirow{2}{*}{pred} & DOFs $1.05$M & $4.20$M & $16.8$M & $67.1$M & $268$M \\
			& & & $h\!=\!9.8\!\times\!10^{-4}$ & $4.9\!\times\!10^{-4}$ & $2.4\!\times\!10^{-4}$ & $1.2\!\times\!10^{-4}$ & $6.1\!\times\!10^{-5}$ \\
			\hline\hline
			\multirow{6}{*}{$2.0004$}
			& \multirow{2}{*}{$\alpha=1$}        & \multirow{2}{*}{$1.5001$}
			& $1.860\!\times\!10^{-4}$ & $6.595\!\times\!10^{-5}$ & $2.336\!\times\!10^{-5}$ & $8.267\!\times\!10^{-6}$ & $2.925\!\times\!10^{-6}$ \\
			& & & $1.49294$ & $1.49581$ & $1.49745$ & $1.49842$ & $1.49900$ \\
			\cline{2-8}
			& \multirow{2}{*}{$\alpha=-\infty$}  & \multirow{2}{*}{$1.5001$}
			& $3.601\!\times\!10^{-4}$ & $1.276\!\times\!10^{-4}$ & $4.517\!\times\!10^{-5}$ & $1.598\!\times\!10^{-5}$ & $5.652\!\times\!10^{-6}$ \\
			& & & $1.49451$ & $1.49698$ & $1.49831$ & $1.49904$ & $1.49945$ \\
			\cline{2-8}
			& \multirow{2}{*}{sym}               & \multirow{2}{*}{$2$}
			& $2.902\!\times\!10^{-6}$ & $7.636\!\times\!10^{-7}$ & $2.003\!\times\!10^{-7}$ & $5.240\!\times\!10^{-8}$ & $1.369\!\times\!10^{-8}$ \\
			& & & $1.91930$ & $1.92606$ & $1.93099$ & $1.93426$ & $1.93616$ \\
			\hline
			\multirow{6}{*}{$4$}
			& \multirow{2}{*}{$\alpha=1$}        & \multirow{2}{*}{$1.75$}
			& $4.490\!\times\!10^{-5}$ & $1.341\!\times\!10^{-5}$ & $4.001\!\times\!10^{-6}$ & $1.193\!\times\!10^{-6}$ & $3.554\!\times\!10^{-7}$ \\
			& & & $1.7402$ & $1.7432$ & $1.7450$ & $1.7461$ & $1.7469$ \\
			\cline{2-8}
			& \multirow{2}{*}{$\alpha=-\infty$}  & \multirow{2}{*}{$1.75$}
			& $8.779\!\times\!10^{-5}$ & $2.617\!\times\!10^{-5}$ & $7.795\!\times\!10^{-6}$ & $2.320\!\times\!10^{-6}$ & $6.905\!\times\!10^{-7}$ \\
			& & & $1.7434$ & $1.7460$ & $1.7474$ & $1.7482$ & $1.7487$ \\
			\cline{2-8}
			& \multirow{2}{*}{sym}               & \multirow{2}{*}{$2$}
			& $2.155\!\times\!10^{-6}$ & $5.407\!\times\!10^{-7}$ & $1.356\!\times\!10^{-7}$ & $3.396\!\times\!10^{-8}$ & $8.500\!\times\!10^{-9}$ \\
			& & & $1.9921$ & $1.9944$ & $1.9959$ & $1.9972$ & $1.9982$ \\
			\hline
			\multirow{6}{*}{$\infty$}
			& \multirow{2}{*}{$\alpha=1$}        & \multirow{2}{*}{$2$}
			& $1.064\!\times\!10^{-5}$ & $2.665\!\times\!10^{-6}$ & $6.668\!\times\!10^{-7}$ & $1.668\!\times\!10^{-7}$ & $4.170\!\times\!10^{-8}$ \\
			& & & $1.9948$ & $1.9974$ & $1.9987$ & $1.9993$ & $1.9997$ \\
			\cline{2-8}
			& \multirow{2}{*}{$\alpha=-\infty$}  & \multirow{2}{*}{$2$}
			& $2.163\!\times\!10^{-5}$ & $5.417\!\times\!10^{-6}$ & $1.355\!\times\!10^{-6}$ & $3.390\!\times\!10^{-7}$ & $8.477\!\times\!10^{-8}$ \\
			& & & $1.9950$ & $1.9975$ & $1.9987$ & $1.9994$ & $1.9997$ \\
			\cline{2-8}
			& \multirow{2}{*}{sym}               & \multirow{2}{*}{$2$}
			& $2.011\!\times\!10^{-6}$ & $5.028\!\times\!10^{-7}$ & $1.257\!\times\!10^{-7}$ & $3.143\!\times\!10^{-8}$ & $7.851\!\times\!10^{-9}$ \\
			& & & $2.0000$ & $2.0000$ & $1.9998$ & $2.0002$ & $2.0010$ \\
			\hline
		\end{tabular}%
	}
\end{table}

\begin{table}[htbp]
	\centering
	\caption{3D convergence histories on the family~\eqref{eq:numerics-family}: relative $L^2$ error (top) and observed rate (bottom). The ``pred'' column is the rate from~\eqref{eq:numerics-predicted-rate}.}
	\label{tab:numerics-3D}
	\resizebox{\textwidth}{!}{%
		\begin{tabular}{c|c|c|ccccc}
			\hline
			\multirow{2}{*}{$p$} & \multirow{2}{*}{Method} & \multirow{2}{*}{pred} & DOFs $4.91$K & $35.9$K & $275$K & $2.15$M & $16.97$M \\
			& & & $h\!=\!6.3\!\times\!10^{-2}$ & $3.1\!\times\!10^{-2}$ & $1.6\!\times\!10^{-2}$ & $7.8\!\times\!10^{-3}$ & $3.9\!\times\!10^{-3}$ \\
			\hline\hline
			\multirow{6}{*}{$2.0004$}
			& \multirow{2}{*}{$\alpha=1$}        & \multirow{2}{*}{$1.5001$}
			& $5.600\!\times\!10^{-2}$ & $2.158\!\times\!10^{-2}$ & $8.072\!\times\!10^{-3}$ & $2.951\!\times\!10^{-3}$ & $1.064\!\times\!10^{-3}$ \\
			& & & $1.33795$ & $1.37559$ & $1.41893$ & $1.45168$ & $1.47215$ \\
			\cline{2-8}
			& \multirow{2}{*}{$\alpha=-\infty$}  & \multirow{2}{*}{$1.5001$}
			& $1.074\!\times\!10^{-1}$ & $4.162\!\times\!10^{-2}$ & $1.553\!\times\!10^{-2}$ & $5.657\!\times\!10^{-3}$ & $2.032\!\times\!10^{-3}$ \\
			& & & $1.28919$ & $1.36775$ & $1.42218$ & $1.45692$ & $1.47695$ \\
			\cline{2-8}
			& \multirow{2}{*}{sym}               & \multirow{2}{*}{$2$}
			& $8.412\!\times\!10^{-3}$ & $2.318\!\times\!10^{-3}$ & $6.291\!\times\!10^{-4}$ & $1.687\!\times\!10^{-4}$ & $4.487\!\times\!10^{-5}$ \\
			& & & $1.83384$ & $1.85935$ & $1.88178$ & $1.89861$ & $1.91093$ \\
			\hline
			\multirow{6}{*}{$4$}
			& \multirow{2}{*}{$\alpha=1$}        & \multirow{2}{*}{$1.75$}
			& $3.859\!\times\!10^{-2}$ & $1.245\!\times\!10^{-2}$ & $3.903\!\times\!10^{-3}$ & $1.197\!\times\!10^{-3}$ & $3.621\!\times\!10^{-4}$ \\
			& & & $1.5395$ & $1.6116$ & $1.6926$ & $1.7259$ & $1.7405$ \\
			\cline{2-8}
			& \multirow{2}{*}{$\alpha=-\infty$}  & \multirow{2}{*}{$1.75$}
			& $7.531\!\times\!10^{-2}$ & $2.448\!\times\!10^{-2}$ & $7.661\!\times\!10^{-3}$ & $2.339\!\times\!10^{-3}$ & $7.043\!\times\!10^{-4}$ \\
			& & & $1.5392$ & $1.6209$ & $1.6762$ & $1.7117$ & $1.7317$ \\
			\cline{2-8}
			& \multirow{2}{*}{sym}               & \multirow{2}{*}{$2$}
			& $9.080\!\times\!10^{-3}$ & $2.346\!\times\!10^{-3}$ & $5.990\!\times\!10^{-4}$ & $1.518\!\times\!10^{-4}$ & $3.831\!\times\!10^{-5}$ \\
			& & & $1.9291$ & $1.9526$ & $1.9694$ & $1.9800$ & $1.9866$ \\
			\hline
			\multirow{6}{*}{$\infty$}
			& \multirow{2}{*}{$\alpha=1$}        & \multirow{2}{*}{$2$}
			& $2.937\!\times\!10^{-2}$ & $8.017\!\times\!10^{-3}$ & $2.118\!\times\!10^{-3}$ & $5.465\!\times\!10^{-4}$ & $1.389\!\times\!10^{-4}$ \\
			& & & $1.8217$ & $1.8734$ & $1.9202$ & $1.9547$ & $1.9758$ \\
			\cline{2-8}
			& \multirow{2}{*}{$\alpha=-\infty$}  & \multirow{2}{*}{$2$}
			& $5.840\!\times\!10^{-2}$ & $1.618\!\times\!10^{-2}$ & $4.289\!\times\!10^{-3}$ & $1.107\!\times\!10^{-3}$ & $2.812\!\times\!10^{-4}$ \\
			& & & $1.7496$ & $1.8521$ & $1.9154$ & $1.9544$ & $1.9763$ \\
			\cline{2-8}
			& \multirow{2}{*}{sym}               & \multirow{2}{*}{$2$}
			& $9.460\!\times\!10^{-3}$ & $2.371\!\times\!10^{-3}$ & $5.933\!\times\!10^{-4}$ & $1.484\!\times\!10^{-4}$ & $3.710\!\times\!10^{-5}$ \\
			& & & $1.9952$ & $1.9965$ & $1.9985$ & $1.9995$ & $1.9998$ \\
			\hline
		\end{tabular}%
	}
\end{table}

\subsection{Sharpness of the threshold \texorpdfstring{$\alpha=1+1/p$}{α=1+1/p}}\label{sec:numerics-critical}

In \Cref{sec:numerics-endpoints} we verified, at the two nonsymmetric endpoints $\alpha=1$ and $\alpha=-\infty$, the predicted rate $r = 2 - 1/p$ across the regularity scale; this confirms that throughout the entire weak-penalty range $\alpha\le 1$, recovering the optimal rate $h^2$ requires the full $W^{2,\infty}$ regularity ($p=\infty$). We now test the critical region $1 < \alpha \le 1.5$, in which \Cref{thm:unified-L2} claims that the optimal rate is recovered if and only if $\alpha \ge 1 + 1/p$:
\begin{align}\label{eq:numerics-critical-predicted}
	r(\alpha,p)\,=\,\min\{2,\,1+\alpha-1/p\}\,=\,\min\{2,\,\alpha+(1-1/p)\}.
\end{align}
As $\alpha$ ranges over $(1,1.5)$, the saturation threshold sweeps the regularity scale from $p=\infty$ down to $p=2$. We take $\alpha\in\{1.1,1.2,1.3,1.4,1.5\}$ and values of $p$ such that $1-1/p\in\{0.5,0.6,0.7,0.8,0.9\}$, so that the saturation boundary $\alpha+(1-1/p)=2$ falls on the anti-diagonal of \Cref{tab:numerics-2D-critical}. From the table we observe entrywise agreement at the finest mesh: sub-optimal cells (sum $<2$, upper-left) match the predicted rate $1+\alpha-1/p$, and saturated cells (sum $\ge 2$) attain $r=2$. This verifies the sharp if-and-only-if threshold of \Cref{thm:unified-L2}: the sub-optimal rate increases linearly in $\alpha$ at each fixed $p$, and increasing the penalty exponent compensates exactly for the regularity deficit $1/p$. To save space, we report only the 2D experiments and do not include $\alpha>3/2$; we have verified that the 3D experiments and the region $\alpha>3/2$ also match the theory exactly.

\begin{table}[htbp]
	\centering
	\caption{2D critical region, nonsymmetric method: observed rate (predicted from~\eqref{eq:numerics-critical-predicted} in parentheses) at the finest mesh ($h\approx 6.1\!\times\!10^{-5}$); $\dagger$ marks $r=2$ saturation. Anti-diagonal: $\alpha+(1-1/p)=2$.}
	\label{tab:numerics-2D-critical}
	\resizebox{\textwidth}{!}{%
		\begin{tabular}{c|ccccc}
			\hline
			\multirow{2}{*}{$\alpha$} & \multicolumn{5}{c}{$1-1/p$} \\
			\cline{2-6}
			& $0.5$ & $0.6$ & $0.7$ & $0.8$ & $0.9$ \\
			\hline
			$1.1$ & $1.5666\ (1.6)$ & $1.6655\ (1.7)$ & $1.7644\ (1.8)$ & $1.8642\ (1.9)$ & $1.9675\ (2)^{\dagger}$ \\
			\hline
			$1.2$ & $1.6649\ (1.7)$ & $1.7636\ (1.8)$ & $1.8631\ (1.9)$ & $1.9660\ (2)^{\dagger}$ & $2.0777\ (2)^{\dagger}$ \\
			\hline
			$1.3$ & $1.7741\ (1.8)$ & $1.8733\ (1.9)$ & $1.9761\ (2)^{\dagger}$ & $2.0871\ (2)^{\dagger}$ & $2.1849\ (2)^{\dagger}$ \\
			\hline
			$1.4$ & $1.8845\ (1.9)$ & $1.9866\ (2)^{\dagger}$ & $2.0951\ (2)^{\dagger}$ & $2.1732\ (2)^{\dagger}$ & $2.0444\ (2)^{\dagger}$ \\
			\hline
			$1.5$ & $1.9950\ (2)^{\dagger}$ & $2.0942\ (2)^{\dagger}$ & $2.1403\ (2)^{\dagger}$ & $2.0154\ (2)^{\dagger}$ & $1.9432\ (2)^{\dagger}$ \\
			\hline
		\end{tabular}%
	}
\end{table}

\paragraph{Conclusion} The numerical experiments confirm the sharpness of \Cref{thm:unified-L2} in both dimensions and across the full penalty scale. In \Cref{sec:numerics-endpoints}, at the two nonsymmetric endpoints $\alpha=1$ and $\alpha=-\infty$, the predicted rate $r=2-1/p$ is realized entrywise across the regularity scale---the half-order loss $h^{3/2}$ at the rough endpoint $p\approx 2$ and the optimal $h^2$ at the smooth endpoint $p=\infty$---while the symmetric reference attains $r=2$ throughout. In \Cref{sec:numerics-critical}, the critical region $1<\alpha\le 1.5$ realizes the predicted rate $\min\{2,\,1+\alpha-1/p\}$ entrywise, with the optimal $r=2$ attained precisely once $\alpha\ge 1+1/p$. These results resolve the long-standing ambiguity: under merely $H^{k+1}$ regularity the half-order loss is essential and intrinsic to the nonsymmetric method, while the optimal rate consistently observed on smooth test problems is explained by their full $W^{k+1,\infty}$ regularity rather than by any limitation of the standard analysis; and the threshold $\alpha\ge 1+1/p$ (or $p=\infty$) is sharp.

\section{Concluding remarks and future directions}\label{sec:conclusion}

On bounded convex polytopes in $\mathbb R^d$, $d\in\{2,3\}$, we proved a unified regularity-dependent $L^2$ error estimate
\begin{equation*}
	\|u-u_h\|_{L^2(\Omega)} \le C\,h^r\,|u|_{W^{k+1,p}(\Omega)},\qquad r=\min\bigl\{k+1,\;k+\max\{1,\alpha\}-1/p\bigr\},
\end{equation*}
uniformly across the entire penalty scale $\gamma=h^{-\alpha}$, $\alpha\in\{-\infty\}\cup\mathbb R$, and the regularity exponent $p\in[2,\infty]$, identifying $\alpha\ge 1+1/p$ or $p=\infty$ as the sharp threshold for the optimal rate $h^{k+1}$. The numerical experiments in \Cref{sec:numerics} confirm the predicted rates entrywise in both $d=2$ and $d=3$ and resolve the open question raised by Burman~\cite[\S 8.1]{MR3022206}: under merely $H^{k+1}$ regularity the half-order loss is essential, while the optimal rate consistently observed in standard tests is explained by the full $W^{k+1,\infty}$ regularity of typical smooth solutions rather than by any limitation of the standard analysis.

The principal restriction of our analysis is the convexity of $\Omega$, which underpins both the dual $H^2$ regularity and the $W^{1,\infty}$ stability of the Ritz projection used in the boundary-strip argument. Removing convexity on general polygonal and polyhedral domains is the natural next step: reentrant corners reduce the dual regularity to $H^{1+s}$ with $s<1$, and recovering the sharp rate will likely require graded meshes near singular points together with weighted-Sobolev arguments tracking the corner components.

\appendix

\section{Proof of \Cref{lem:disc-lift}}\label{app:lem:disc-lift}
\begin{proof}[Proof of \Cref{lem:disc-lift}]
	Let $\mathcal N_h^\Gamma$ denote the set of Lagrange nodes whose geometric position lies on $\Gamma$; this is precisely the set of degrees of freedom of $M_h=V_h|_\Gamma$.  In addition to the vertex nodes on $\Gamma$, for $k\ge 2$ it contains edge-interior Lagrange nodes lying on $\Gamma$, and in three dimensions for $k\ge 3$ also face-interior Lagrange nodes lying on $\Gamma$.  Let $\varphi_a$ denote the standard nodal Lagrange basis function of $V_h$ at $a\in\mathcal N_h^\Gamma$.  Define
	\begin{equation*}
		E_h\mu_h:=\sum_{a\in\mathcal N_h^\Gamma}\mu_h(a)\varphi_a. 
	\end{equation*}
	Then $(E_h\mu_h)|_\Gamma=\mu_h$ (since $M_h$ has $\mathcal N_h^\Gamma$ as its degrees of freedom) and $\operatorname{supp}(E_h\mu_h)\subset\overline{\Omega_h^\Gamma}$ (for each $a\in\mathcal N_h^\Gamma$, $\varphi_a$ is supported in the union of closed mesh elements containing $a$; since $a\in\Gamma$, every such element $T$ satisfies $\overline T\cap\Gamma\ne\emptyset$, hence $T\in\mathcal T_h^\Gamma$, so $\overline T\subset\overline{\Omega_h^\Gamma}$ and $\operatorname{supp}(E_h\mu_h)\subset\overline{\Omega_h^\Gamma}$; since $\overline{\Omega_h^\Gamma}\setminus\Omega_h^\Gamma$ has Lebesgue measure zero, $E_h\mu_h=0$ a.e.\ in $\Omega\setminus\Omega_h^\Gamma$), giving~\eqref{eq:disc-lift-trace}.
	
	Finite-dimensional norm equivalence on the reference element together with shape-regular scaling gives
	\begin{equation}\label{eq:disc-lift-patch}
		\|\nabla E_h\mu_h\|_{L^2(T)}^2
		\le C h_T^{d-2}\sum_{a\in\mathcal N_h^\Gamma\cap\overline T}|\mu_h(a)|^2,\qquad T\in\mathcal T_h^\Gamma. 
	\end{equation}
	For each boundary face $E\in\mathcal E_h^\Gamma$, the analogous norm equivalence on the reference face together with shape-regular scaling gives
	\begin{equation*}
		h^{d-1}\sum_{a\in\mathcal N_h^\Gamma\cap\overline E}|\mu_h(a)|^2\le C\|\mu_h\|_{L^2(E)}^2. 
	\end{equation*}
	Summing over $E\in\mathcal E_h^\Gamma$ and using that every node in $\mathcal N_h^\Gamma$ lies on at least one boundary face yields
	\begin{equation}\label{eq:disc-lift-bnd-norm-eq}
		h^{d-1}\sum_{a\in\mathcal N_h^\Gamma}|\mu_h(a)|^2
		\le C\|\mu_h\|_{L^2(\Gamma)}^2. 
	\end{equation}
	Summing~\eqref{eq:disc-lift-patch} over $T\in\mathcal T_h^\Gamma$, using the quasi-uniformity $h_T\simeq h$ and the fact that each Lagrange node belongs to a uniformly bounded number of elements (by shape regularity), and combining with~\eqref{eq:disc-lift-bnd-norm-eq} yields
	\begin{equation*}
		\|\nabla E_h\mu_h\|_{L^2(\Omega)}^2
		\le Ch^{d-2}\sum_{a\in\mathcal N_h^\Gamma}|\mu_h(a)|^2
		\le Ch^{-1}\|\mu_h\|_{L^2(\Gamma)}^2,
	\end{equation*}
	which gives~\eqref{eq:disc-lift-W12} on taking square roots.  Finally, the inverse trace~\eqref{eq:grad-inverse-trace} applied face by face gives
	\begin{equation*}
		\|\partial_n E_h\mu_h\|_{L^2(\Gamma)}
		\le Ch^{-1/2}\|\nabla E_h\mu_h\|_{L^2(\Omega_h^\Gamma)}
		\le Ch^{-1}\|\mu_h\|_{L^2(\Gamma)}.
	\end{equation*}
\end{proof}

\section{Proof of \Cref{ritz_pro_L2err}}\label{app:ritz_pro_L2err}
\begin{proof}[Proof of \Cref{ritz_pro_L2err}]
	Let $e_I:=v-I_hv$.  We first record the two boundary interpolation estimates
	\begin{equation}\label{eq:boundary-interp-basic}
		\|e_I\|_{L^2(\Gamma)}\le Ch^{k+1/2}|v|_{H^{k+1}(\Omega)},
		\qquad
		\|\partial_ne_I\|_{L^2(\Gamma)}\le Ch^{k-1/2}|v|_{H^{k+1}(\Omega)}. 
	\end{equation}
	For a boundary face $E\subset\partial T$, the scaled trace inequality~\eqref{eq:scaled-trace} and the interpolation estimate~\eqref{eq:H-interp-volume} give
	\begin{equation*}
		\|e_I\|_{L^2(E)}^2
		\le C\bigl(h_T^{-1}\|e_I\|_{L^2(T)}^2+h_T\|\nabla e_I\|_{L^2(T)}^2\bigr)
		\le C h_T^{2k+1}|v|_{H^{k+1}(T)}^2,
	\end{equation*}
	and summing over boundary faces proves the first bound.  Applying~\eqref{eq:scaled-trace} to $\nabla e_I$ and using~\eqref{eq:H-interp-volume} with $m=1,2$ gives
	\begin{equation*}
		\|\partial_ne_I\|_{L^2(E)}^2
		\le C\bigl(h_T^{-1}\|\nabla e_I\|_{L^2(T)}^2+h_T\|D^2 e_I\|_{L^2(T)}^2\bigr)
		\le Ch_T^{2k-1}|v|_{H^{k+1}(T)}^2,
	\end{equation*}
	proving the second bound.
	
	\smallskip\noindent\emph{Bound \eqref{eq:Ritz-H-boundary}.}
	By construction $(\mathcal R_hv)|_\Gamma=\Pi_h^\Gamma(v|_\Gamma)$, and $(I_hv)|_\Gamma\in M_h$.  The $L^2(\Gamma)$ best-approximation property of $\Pi_h^\Gamma$ therefore yields
	\begin{equation*}
		\|v-\mathcal R_hv\|_{L^2(\Gamma)}
		=\|v|_\Gamma-\Pi_h^\Gamma(v|_\Gamma)\|_{L^2(\Gamma)}
		\le\|e_I\|_{L^2(\Gamma)}\le Ch^{k+1/2}|v|_{H^{k+1}(\Omega)}. 
	\end{equation*}
	
	\smallskip\noindent\emph{Bound~\eqref{eq:Ritz-H-global}.}
	Set $\mu_h:=\Pi_h^\Gamma(v|_\Gamma)-(I_hv)|_\Gamma$.  Since $(I_hv)|_\Gamma\in M_h$, $\mu_h=\Pi_h^\Gamma(e_I|_\Gamma)$, and the $L^2(\Gamma)$ stability~\eqref{eq:Pi-L2-stab} gives $\|\mu_h\|_{L^2(\Gamma)}\le\|e_I\|_{L^2(\Gamma)}$.  Define $J_hv:=I_hv+E_h\mu_h\in V_h$, with $E_h$ the lifting of \Cref{lem:disc-lift}; then $(J_hv)|_\Gamma=(I_hv)|_\Gamma+\mu_h=(\mathcal R_hv)|_\Gamma$, so $\mathcal R_hv-J_hv\in V_h^0$.  Ritz orthogonality in~\eqref{eq:nh-ritz-def} applied to $w_h=\mathcal R_hv-J_hv$ yields the best-approximation
	\begin{equation*}
		\|\nabla(v-\mathcal R_hv)\|_{L^2(\Omega)}\le\|\nabla(v-J_hv)\|_{L^2(\Omega)}. 
	\end{equation*}
	Combining~\eqref{eq:H-interp-volume},~\eqref{eq:disc-lift-W12}, and~\eqref{eq:boundary-interp-basic} gives
	\begin{align*}
		\|\nabla(v-J_hv)\|_{L^2(\Omega)}
		&\le\|\nabla(v-I_hv)\|_{L^2(\Omega)}+\|\nabla E_h\mu_h\|_{L^2(\Omega)} \\
		&\le Ch^k|v|_{H^{k+1}(\Omega)}+Ch^{-1/2}\|e_I\|_{L^2(\Gamma)}
		\le Ch^k|v|_{H^{k+1}(\Omega)}. 
	\end{align*}
	
	\smallskip\noindent\emph{Bound~\eqref{eq:Ritz-H-normal}.}
	The triangle inequality yields
	\begin{equation*}
		\|\partial_n(v-\mathcal R_hv)\|_{L^2(\Gamma)}
		\le\|\partial_ne_I\|_{L^2(\Gamma)}+\|\partial_n(I_hv-\mathcal R_hv)\|_{L^2(\Gamma)}. 
	\end{equation*}
	The first term is bounded by~\eqref{eq:boundary-interp-basic}.  For the second, $I_hv-\mathcal R_hv\in V_h$, so the inverse trace~\eqref{eq:grad-inverse-trace}, the triangle inequality, \eqref{eq:H-interp-volume}, and the already-proved bound~\eqref{eq:Ritz-H-global} give
	\begin{equation*}
		\|\partial_n(I_hv-\mathcal R_hv)\|_{L^2(\Gamma)}
		\le Ch^{-1/2}\|\nabla(I_hv-\mathcal R_hv)\|_{L^2(\Omega)}
		\le Ch^{k-1/2}|v|_{H^{k+1}(\Omega)}. 
	\end{equation*}
\end{proof}

\section{Proof of \Cref{thm:weak-infsup}}\label{app:weak-infsup}

The proof is inspired by Burman's inf--sup philosophy~\cite{MR3022206}, but the local construction differs. Since Burman's analysis is presented in two dimensions and does not directly give a 3D proof, for completeness we provide a rigorous proof uniformly for both $d=2$ and $d=3$.

Throughout this appendix we use the penalty-free part of the bilinear form and the reference weak-penalty norm
\begin{equation}\label{eq:app-a0-norm}
	\begin{split}
		a_0(v,w)&:=(\nabla v,\nabla w)_\Omega-\langle\partial_n v,w\rangle_\Gamma+\langle v,\partial_n w\rangle_\Gamma,\\
		\vertiii{v}_{1,h}^2&:=\|\nabla v\|_{L^2(\Omega)}^2+h^{-1}\|v\|_{L^2(\Gamma)}^2 .
	\end{split}
\end{equation}
For $\alpha<1$, $\vertiii{\cdot}_{1,h}$ is the norm $\vertiii{\cdot}_h$ in~\eqref{eq:reference-nitsche-norm}. Constants below may depend on the polynomial degree, the dimension, the domain, and the mesh-regularity constants, but not on $h$ or on the finite value of $\alpha<1$.

For every boundary face $E \in \mathcal{E}_h^\Gamma$, let $T_E$ be the unique parent simplex with $E \subset \partial T_E$. The unique vertex of $T_E$ that does not lie on the closed face $\overline{E}$, denoted by $z_E$ and called the opposite vertex of $E$.
Let $\mathfrak S$ denote the lower-dimensional boundary skeleton: the set of polygonal vertices if $d=2$, and the union of polyhedral edges and vertices if $d=3$. Split the boundary faces into
\begin{equation}\label{eq:app-std-exc-def}
	\mathcal S_h:=\{E\in\mathcal E_h^\Gamma:\text{ the opposite vertex }z_E\in\Omega\setminus\Gamma\},
	\qquad
	\mathcal X_h:=\mathcal E_h^\Gamma\setminus\mathcal S_h .
\end{equation}
Faces in $\mathcal S_h$ are called standard; faces in $\mathcal X_h$ are called exceptional. Set
\[
\Gamma_S:=\bigcup_{E\in\mathcal S_h}E,
\qquad
\Gamma_X:=\bigcup_{E\in\mathcal X_h}E .
\]
Standard faces are controlled by the normal derivative of the affine hat attached to their opposite interior vertex.  Exceptional faces are the corner and edge configurations where no such
interior opposite vertex is available. These faces are localized in an
\(O(h)\)-neighborhood of the boundary skeleton. For each exceptional face
we construct a uniformly bounded macro-patch of diameter \(O(h)\), connect
it to a neighboring standard face, and transfer the trace on the exceptional
face to the corresponding interior anchor. The resulting macro-patches have
uniformly bounded overlap, and each anchor is used only a uniformly bounded
number of times.

The following elementary geometric lemma is the rigorous form of the corner/edge macro-patch construction. It is the only point where the threshold $h_0$ enters.

\begin{lemma}[corner/edge macro-patch]\label{lem:app-exc-face-patch}
	There is a number $h_0>0$, depending only on $\Omega$ and on the shape-regularity and quasi-uniformity constants, with the following property. On every mesh with $h<h_0$, each exceptional face $E\in\mathcal X_h$ admits a face-connected element macro-patch $\omega_E$ and an interior anchor $A_E$. The anchor is the opposite interior vertex of at least one neighboring standard face $F_E\in\mathcal S_h$. They may be chosen so that
	\begin{equation}\label{eq:app-patch-size}
		\operatorname{diam}\omega_E\le Ch,
		\qquad
		\#\{T\in\mathcal T_h:T\subset\omega_E\}\le C,
		\qquad
		A_E\in\omega_E,
	\end{equation}
	that the patches have uniformly bounded overlap, and that the anchors have uniformly bounded multiplicity:
	\begin{equation}\label{eq:app-anchor-multiplicity}
		\sup_{z\in\{z_F:F\in\mathcal S_h\}}\#\{E\in\mathcal X_h:A_E=z\}\le C .
	\end{equation}
	Moreover, for every $v_h\in V_h$,
	\begin{equation}\label{eq:app-exc-face-control}
		h^{-1}\|v_h\|_{L^2(E)}^2
		\le C\left(
		\|\nabla v_h\|_{L^2(\omega_E)}^2
		+h^{d-2}|v_h(A_E)|^2
		\right).
	\end{equation}
\end{lemma}

\begin{proof}
	Let $\mathfrak F$ be the  set of flat boundary facets of $\Omega$ (edges if $d=2$, polygonal faces if $d=3$) For $F\in\mathfrak F$, write $\partial_\Gamma F$ for its relative boundary in $\Gamma$. Thus $\bigcup_{F\in\mathfrak F}\partial_\Gamma F=\mathfrak S$.
	
	We first show that exceptional faces are confined to a boundary strip of thickness $O(h)$ around the skeleton. There is a constant $C_\Omega$ such that, whenever a boundary mesh face $E\subset F$ is exceptional,
	\begin{equation}\label{eq:app-exc-near-skeleton}
		\operatorname{dist}(E,\partial_\Gamma F)\le C_\Omega h .
	\end{equation}
	Indeed, let $z_E$ be the opposite vertex of $T_E$. Since $E$ is exceptional, $z_E\in\Gamma$. Since all vertices of $E$ lie in the single flat facet $F$, the point $z_E$ cannot also lie in $F$ unless $T_E$ is degenerate. Hence $z_E\in\Gamma\setminus F$. Because $\operatorname{diam}T_E\le h$, $\operatorname{dist}(E,\Gamma\setminus F)\le h$. On a fixed convex polytope the finite set of facets has a uniform angle separation; equivalently, for each facet $F$ there is $c_F>0$ such that
	\begin{equation}\label{eq:app-angle-distance}
		\operatorname{dist}(x,\Gamma\setminus F)
		\ge c_F\operatorname{dist}(x,\partial_\Gamma F)
		\qquad x\in F .
	\end{equation}
	Taking the minimum over the finitely many facets proves~\eqref{eq:app-exc-near-skeleton}. In particular, every boundary mesh face $G\subset F$ with $\operatorname{dist}(G,\partial_\Gamma F)>C_\Omega h$ is standard.

	Fix a exceptional face $E\subset F$. Choose $x_E\in E$. Since $F$ is a fixed convex facet and $h_0$ is chosen below its inradius scale, there are constants $C_1,C_2$, depending only on $\Omega$, and a point $y_E\in F$ such that
	\[
	|x_E-y_E|\le C_1h,\qquad \operatorname{dist}(y_E,\partial_\Gamma F)\ge (C_\Omega+2C_2)h .
	\]
	Let $F_E\subset F$ be the boundary mesh face containing $y_E$ in its closure. Since the boundary mesh is quasi-uniform and $\operatorname{diam}F_E\le C_2h$ after enlarging $C_2$ if necessary,
	\begin{equation}\label{eq:app-std-face-away-skeleton}
		\operatorname{dist}(F_E,\partial_\Gamma F)>C_\Omega h .
	\end{equation}
	By the preceding paragraph, $F_E$ is standard. Let $A_E:=z_{F_E}\in\Omega\setminus\Gamma$ be its opposite interior vertex.

	The boundary mesh induced on each fixed facet is uniformly shape regular and quasi-uniform. Hence $E$ and $F_E$ can be connected inside $F\cap B(E,Ch)$ by a chain of at most $C$ boundary mesh faces. Let $\omega_E$ initially be the union of the parent elements of this chain together with the elements needed to include the anchor $A_E$. If two consecutive parent elements in the boundary chain meet only through a lower-dimensional boundary entity {(e.g., a vertex in 2D, vertexes and an edge in 3D)}, we enlarge the patch by the uniformly bounded element star of that entity. Shape regularity gives a uniformly positive lower bound on the corresponding element angles, so this star contains only $O(1)$ elements. After this enlargement $\omega_E$ is face-connected, i.e., any two elements in the patch are linked by a chain of elements sharing full $(d-1)$-dimensional faces. This gives the diameter and cardinality bounds in~\eqref{eq:app-patch-size}. The bounded-overlap property and the anchor multiplicity bound~\eqref{eq:app-anchor-multiplicity} follow from quasi-uniformity, since every construction is contained in a boundary ball of radius $Ch$ and contains only $O(1)$ elements.
	
	Scale $\omega_E$ to diameter one. On the scaled patch the relevant finite element space has uniformly bounded dimension and uniformly shape-regular geometry. The reference-patch inequality
	\[
	\|\widehat v_h\|_{L^2(\widehat E)}^2
	\le C\left(
	\|\widehat\nabla\widehat v_h\|_{L^2(\widehat\omega_E)}^2
	+|\widehat v_h(\widehat A_E)|^2
	\right)
	\]
	follows by finite-dimensional norm equivalence. If the right-hand side is zero, then $\widehat v_h$ is constant on each element; because the patch is face-connected and functions in $V_h$ have matching traces across common element faces, these elementwise constants are equal throughout the patch. The condition $\widehat v_h(\widehat A_E)=0$ then forces the common constant to be zero. Uniformity of the constant follows from compactness of the uniformly shape-regular patches with bounded cardinality and fixed polynomial degree. Scaling back to diameter $h$ gives~\eqref{eq:app-exc-face-control}.
\end{proof}

We now organize the standard boundary faces according to their opposite interior vertices. Let {$\{z_j\}_{j\in \mathcal J}$, where $\mathcal{J}$ is a finite index set,} enumerate the distinct opposite vertices $\{z_E:E\in\mathcal S_h\}$ and define
\begin{equation}\label{eq:app-boundary-patch-def}
	\mathcal E_j:=\{E\in\mathcal S_h:\ z_E=z_j\},
	\qquad
	\Gamma_j:=\bigcup_{E\in\mathcal E_j}E,
\end{equation}
\begin{equation}\label{eq:app-volume-star-def}
	\omega_j:=\operatorname{int}\bigcup_{\{T\in\mathcal T_h:\ z_j\in T\}}\overline T .
\end{equation}
The sets $\Gamma_j$ cover $\Gamma_S$ and are pairwise disjoint up to sets of surface measure zero. By shape regularity and quasi-uniformity, the number of elements in each vertex star is uniformly bounded, the stars $\omega_j$ have uniformly bounded overlap, and
\begin{equation}\label{eq:app-patch-scaling}
	\operatorname{diam}\Gamma_j\le Ch,\qquad
	c h^{d-1}\le |\Gamma_j|\le C h^{d-1},\qquad
	|\omega_j|\le Ch^d .
\end{equation}

For each $z_j$, let $\psi_j$ be the global continuous piecewise affine nodal hat function on $\mathcal T_h$ satisfying $\psi_j(z_j)=1$ and $\psi_j(b)=0$ at every other mesh vertex $b$. Since $k\ge1$, this affine hat function belongs to $V_h$.

\begin{lemma}\label{lem:app-hat-geometry}
	For each standard-face vertex $z_j$,
	\begin{equation}\label{eq:app-hat-trace}
		\psi_j|_\Gamma=0,
	\end{equation}
	and, on every face $E\subset\Gamma_j$,
	\begin{equation}\label{eq:app-eta-bounds}
		c h^{-1}\le (-\partial_n\psi_j)|_E\le C h^{-1}.
	\end{equation}
	If $i\ne j$, then $\partial_n\psi_i=0$ on $\Gamma_j$, and every $\partial_n\psi_i$ also vanishes on $\Gamma_X$. Moreover,
	\begin{equation}\label{eq:app-hat-energy}
		\|\nabla\psi_j\|_{L^2(\omega_j)}^2\le C h^{d-2},
		\qquad
		c h^{d-2}\le M_j:=\int_{\Gamma_j}(-\partial_n\psi_j)\,ds\le C h^{d-2}.
	\end{equation}
\end{lemma}

\begin{proof}
	Because $z_j$ is an interior vertex, it is not a vertex of any boundary face. Hence the affine function $\psi_j$ is zero at all vertices of every boundary face and therefore vanishes identically on $\Gamma$. If $E\subset\Gamma_j$, then $z_j$ is the vertex of $T_E$ opposite to $E$, so $\psi_j|_{T_E}$ is the affine barycentric coordinate associated with this opposite vertex. It vanishes on $E$ and increases in the inward normal direction. Thus $-\partial_n\psi_j$ is the reciprocal of the altitude from $z_j$ to $E$, and shape regularity plus quasi-uniformity give~\eqref{eq:app-eta-bounds}. If $E\subset\Gamma_j$ and $i\ne j$, then $z_i$ is not a vertex of $T_E$, so $\psi_i$ vanishes on $T_E$. If $E\in\mathcal X_h$, all vertices of $T_E$ lie on $\Gamma$, so no standard-face vertex belongs to $T_E$ and every $\psi_i$ vanishes on $T_E$. The remaining estimates follow from $|\nabla\psi_j|\le Ch^{-1}$ on $\omega_j$, from~\eqref{eq:app-patch-scaling}, and from~\eqref{eq:app-eta-bounds}.
\end{proof}

\begin{lemma}\label{lem:app-vertex-star-boundary}
	Let $E\in\mathcal E_j$ and let $T_E$ be its parent element. For all $v_h\in V_h$,
	\begin{equation}\label{eq:app-face-trace-with-vertex}
		h^{-1}\|v_h\|_{L^2(E)}^2
		\le C\left(\|\nabla v_h\|_{L^2(T_E)}^2+h^{d-2}|v_h(z_j)|^2\right),
	\end{equation}
	and
	\begin{equation}\label{eq:app-face-trace-zero-at-vertex}
		h^{-1}\|v_h-v_h(z_j)\|_{L^2(E)}^2
		\le C\|\nabla v_h\|_{L^2(T_E)}^2.
	\end{equation}
\end{lemma}

\begin{proof}
	The estimates are finite-dimensional scaling on a simplex. On a reference simplex with opposite face $\widehat E$ and opposite vertex $\widehat z$, the quantity
	\[
	\|\widehat\nabla\widehat v\|_{L^2(\widehat T)}^2+|\widehat v(\widehat z)|^2
	\]
	is a squared norm on $\mathcal P_k(\widehat T)$; {
		on the subspace $\{\hat{v} \in \mathcal{P}_k(\hat{T}) : \hat{v}(\hat{z}) = 0\}$, the seminorm $\|\hat{\nabla}\hat{v}\|_{L^2(\hat{T})}$ is a norm, this gives
		\begin{align*}
			\|\hat{v}_h - \hat{v}_h(\hat{z})\|_{L^2(\hat{E})}^2 \le C\|\nabla \hat{v}_h\|_{L^2(\hat{T})}^2.
		\end{align*}
	} Pulling these two reference estimates back to $T_E$ gives~\eqref{eq:app-face-trace-with-vertex} and~\eqref{eq:app-face-trace-zero-at-vertex}.
\end{proof}

\begin{lemma}\label{lem:app-weighted-average-control}
	For $v_h\in V_h$, define
	\begin{equation}\label{eq:app-weighted-averages}
		M_j:=\int_{\Gamma_j}(-\partial_n\psi_j)\,ds,\qquad
		\bar v_j:=M_j^{-1}\int_{\Gamma_j}(-\partial_n\psi_j) v_h\,ds,\qquad
		\Sigma(v_h):=\sum_{j\in {\mathcal{J}}}M_j\bar v_j^2 .
	\end{equation}
	Then
	\begin{equation}\label{eq:app-std-boundary-control}
		h^{-1}\|v_h\|_{L^2(\Gamma_S)}^2
		\le C\left(\|\nabla v_h\|_{L^2(\Omega)}^2+\Sigma(v_h)\right),
	\end{equation}
	and therefore
	\begin{equation}\label{eq:app-full-boundary-control}
		h^{-1}\|v_h\|_{L^2(\Gamma)}^2
		\le C\left(\|\nabla v_h\|_{L^2(\Omega)}^2+\Sigma(v_h)\right).
	\end{equation}
\end{lemma}

\begin{proof}
	{First, we control the standard part.} By Jensen's inequality and~\eqref{eq:app-face-trace-zero-at-vertex},
	\begin{align}
		M_j|\bar v_j-v_h(z_j)|^2
		&\le \int_{\Gamma_j}(-\partial_n\psi_j)|v_h-v_h(z_j)|^2\,ds \notag\\
		&\le C\sum_{E\in\mathcal E_j}h^{-1}\|v_h-v_h(z_j)\|_{L^2(E)}^2
		\le C\|\nabla v_h\|_{L^2(\omega_j)}^2 . \label{eq:app-average-controls-vertex}
	\end{align}
	Since $M_j\ge ch^{d-2}$, this gives
	\begin{equation}\label{eq:app-vertex-by-average}
		h^{d-2}|v_h(z_j)|^2
		\le C\left(M_j\bar v_j^2+\|\nabla v_h\|_{L^2(\omega_j)}^2\right).
	\end{equation}
	Summing~\eqref{eq:app-face-trace-with-vertex} over $E\in\mathcal E_j$ and using~\eqref{eq:app-vertex-by-average} gives
	\[
	h^{-1}\|v_h\|_{L^2(\Gamma_j)}^2
	\le C\left(\|\nabla v_h\|_{L^2(\omega_j)}^2+M_j\bar v_j^2\right).
	\]
	Summing over $j$ and using bounded overlap proves~\eqref{eq:app-std-boundary-control}.
	
	It remains to control the exceptional part. For each $E\in\mathcal X_h$, the macro-patch lemma gives an anchor $A_E=z_{j(E)}$ for some $j(E)\in {\mathcal{J}}$. Combining~\eqref{eq:app-exc-face-control} with~\eqref{eq:app-vertex-by-average},
	\[
	h^{-1}\|v_h\|_{L^2(E)}^2
	\le C\left(
	\|\nabla v_h\|_{L^2(\omega_E)}^2
	+M_{j(E)}\bar v_{j(E)}^2
	+\|\nabla v_h\|_{L^2(\omega_{j(E)})}^2
	\right).
	\]
	Summing over $E\in\mathcal X_h$, {using the bounded overlap of the macro-patches $\omega_E$ and the vertex stars $\omega_{j(E)}$}, and the anchor multiplicity bound~\eqref{eq:app-anchor-multiplicity}, gives
	\begin{equation}\label{eq:app-exc-boundary-by-anchors}
		h^{-1}\|v_h\|_{L^2(\Gamma_X)}^2
		\le C\left(\|\nabla v_h\|_{L^2(\Omega)}^2+\Sigma(v_h)\right).
	\end{equation}
	Combining~\eqref{eq:app-std-boundary-control} and~\eqref{eq:app-exc-boundary-by-anchors} proves~\eqref{eq:app-full-boundary-control}.
\end{proof}

\begin{lemma}\label{lem:app-flux-test-function}
	For each $v_h\in V_h$, define
	\begin{equation}\label{eq:app-phi-def-new}
		\phi_h:=-\sum_{j\in {\mathcal{J}}}\bar v_j\psi_j,
	\end{equation}
	where $\bar v_j$ is given by~\eqref{eq:app-weighted-averages}. Then $\phi_h\in V_h$, $\phi_h|_\Gamma=0$, and
	\begin{equation}\label{eq:app-flux-main-new}
		\langle v_h,\partial_n\phi_h\rangle_\Gamma=\Sigma(v_h),
	\end{equation}
	\begin{equation}\label{eq:app-phi-stability-new}
		\vertiii{\phi_h}_{1,h}^2=\|\nabla\phi_h\|_{L^2(\Omega)}^2\le C \Sigma(v_h).
	\end{equation}
\end{lemma}

\begin{proof}
	{
		The property $\phi_h|_{\Gamma} = 0$ follows directly from \Cref{lem:app-hat-geometry}. Consequently, the boundary contribution to the $\|\phi_h\|_{1,h}$ vanishes, yielding $\vertiii{\phi_h}_{1,h}^2 = \|\nabla \phi_h\|_{L^2(\Omega)}^2$.
		
		On standard faces $\Gamma_S = \bigcup_{j \in \mathcal{J}} \Gamma_j$, $\psi_i|_{\Gamma_j} = 0$ for all $i \neq j$, hence $\partial_n \phi_h|_{\Gamma_j} = -\bar{v}_j \partial_n \psi_j$; on $\Gamma_X$, all $\psi_j$ ($j \in \mathcal{J}$) vanish in a neighborhood, so $\partial_n \phi_h|_{\Gamma_X} = 0$.
	}
	Hence
	\[
	\langle v_h,\partial_n\phi_h\rangle_\Gamma
	=\sum_{j\in {\mathcal{J}}}\bar v_j\int_{\Gamma_j}(-\partial_n\psi_j)v_h\,ds
	=\Sigma(v_h).
	\]
	For stability, at most $d+1$ affine hats are nonzero on any element, and the vertex stars have bounded overlap. Hence, using~\eqref{eq:app-hat-energy},
	\[
	{\vertiii{\phi_h}_{1,h}^2}=\|\nabla\phi_h\|_{L^2(\Omega)}^2
	\le C\sum_{j\in {\mathcal{J}}}\bar v_j^2\|\nabla\psi_j\|_{L^2(\omega_j)}^2
	\le C\sum_{j\in {\mathcal{J}}}M_j\bar v_j^2=C \Sigma(v_h).
	\]
\end{proof}

\begin{theorem}\label{thm:app-penalty-free-infsup}
	On every mesh with $h<h_0$, there exists $c_0>0$, independent of $h$, such that, for all $v_h\in V_h$,
	\begin{equation}\label{eq:app-penalty-free-infsup}
		c_0\vertiii{v_h}_{1,h}
		\le
		\sup_{0\ne w_h\in V_h}\frac{|a_0(v_h,w_h)|}{\vertiii{w_h}_{1,h}} .
	\end{equation}
\end{theorem}

\begin{proof}
	Let $\phi_h$ be defined by~\eqref{eq:app-phi-def-new}. Since $\phi_h|_\Gamma=0$,
	\begin{equation}\label{eq:app-a0-v-phi-new}
		a_0(v_h,\phi_h)=(\nabla v_h,\nabla\phi_h)_\Omega+\langle v_h,\partial_n\phi_h\rangle_\Gamma
		=(\nabla v_h,\nabla\phi_h)_\Omega+\Sigma(v_h).
	\end{equation}
	Choose $w_h:=v_h+\delta\phi_h$, where $\delta>0$ is fixed independently of $h$. The diagonal cancellation gives $a_0(v_h,v_h)=\|\nabla v_h\|_{L^2(\Omega)}^2$, and therefore
	\begin{equation}\label{eq:app-a0-v-w-new}
		a_0(v_h,w_h)
		=\|\nabla v_h\|_{L^2(\Omega)}^2
		+\delta \Sigma(v_h)+\delta(\nabla v_h,\nabla\phi_h)_\Omega .
	\end{equation}
	By~\eqref{eq:app-phi-stability-new} and Young's inequality,
	\[
	\delta |(\nabla v_h,\nabla\phi_h)_\Omega|
	\le \frac12\|\nabla v_h\|_{L^2(\Omega)}^2+C\delta^2\Sigma(v_h).
	\]
	Taking $0<\delta\le\delta_0$ small enough gives
	\begin{equation}\label{eq:app-a0-lower-new}
		a_0(v_h,w_h)
		\ge c\left(\|\nabla v_h\|_{L^2(\Omega)}^2+\Sigma(v_h)\right).
	\end{equation}
	Together with~\eqref{eq:app-full-boundary-control}, this yields
	\begin{equation}\label{eq:app-control-full-norm-new}
		a_0(v_h,w_h)\ge c\vertiii{v_h}_{1,h}^2 .
	\end{equation}
	It remains to bound the norm of the chosen test function. {By Cauchy's inequality for the terms $(-\partial_n\psi_j)^{1/2}$ and $(-\partial_n\psi_j)^{1/2}v_h$, yields}
	\begin{equation}\label{eq:app-S-upper}
		\Sigma(v_h)=\sum_{j\in {\mathcal{J}}}\frac{\left(\int_{\Gamma_j}(-\partial_n\psi_j)v_h\,ds\right)^2}{M_j}
		\le \sum_{j\in {\mathcal{J}}}\int_{\Gamma_j}(-\partial_n\psi_j)v_h^2\,ds
		\le C h^{-1}\|v_h\|_{L^2(\Gamma)}^2.
	\end{equation}
	Therefore, by~\eqref{eq:app-phi-stability-new},
	\begin{equation}\label{eq:app-w-bound-new}
		\vertiii{w_h}_{1,h}
		\le \vertiii{v_h}_{1,h}+\delta\vertiii{\phi_h}_{1,h}
		\le C\vertiii{v_h}_{1,h}.
	\end{equation}
	Combining~\eqref{eq:app-control-full-norm-new} and~\eqref{eq:app-w-bound-new} proves~\eqref{eq:app-penalty-free-infsup}.
\end{proof}

\begin{proof}[Proof of \Cref{thm:weak-infsup}]
	For $\alpha=-\infty$, the result is \Cref{thm:app-penalty-free-infsup}. Let now $\alpha<1$ be finite. Then
	\[
	a(v,w)=a_0(v,w)+h^{-\alpha}\langle v,w\rangle_\Gamma,
	\qquad
	\vertiii{\cdot}_h=\vertiii{\cdot}_{1,h}.
	\]
	Use the same test function $w_h=v_h+\delta\phi_h$ as in the penalty-free proof. Since $\phi_h|_\Gamma=0$, one has $w_h|_\Gamma=v_h|_\Gamma$, and hence
	\[
	h^{-\alpha}\langle v_h,w_h\rangle_\Gamma=h^{-\alpha}\|v_h\|_{L^2(\Gamma)}^2\ge0.
	\]
	Consequently,
	\[
	a(v_h,w_h)
	=a_0(v_h,w_h)+h^{-\alpha}\|v_h\|_{L^2(\Gamma)}^2
	\ge c\vertiii{v_h}_{1,h}^2.
	\]
	The bound $\vertiii{w_h}_{1,h}\le C\vertiii{v_h}_{1,h}$ is unchanged. Therefore
	\[
	\sup_{0\ne w_h\in V_h}\frac{|a(v_h,w_h)|}{\vertiii{w_h}_h}
	\ge c\vertiii{v_h}_h .
	\]
	The constant is independent of the finite value of $\alpha<1$.
\end{proof}

\bibliographystyle{siamplain}
\bibliography{references}

\end{document}